\documentclass[preprint,3p]{elsarticle}

\usepackage{amssymb}
\usepackage{amsmath}
\usepackage{algorithm}
\usepackage{algorithmic}
\usepackage{color}
\usepackage{booktabs}
\newcommand{\bs}[1]{\boldsymbol{#1}}

\newcommand{\figref}[1]{{Fig. \ref{#1}}}

\begin{document}

\begin{frontmatter}



\title{Point Jacobi-type preconditioning and parameter tuning for Calderon-preconditioned Burton--Miller method in transmission problems}


\author[label1]{Keigo Tomoyasu} %% Author name
\author[label1]{Hiroshi Isakari\corref{cor1}} %% Author name
\ead{isakari@sd.keio.ac.jp}
\cortext[cor1]{Corresponding author}

\affiliation[label1]{organization={Faculty of Science and Technology, Keio University},%Department and Organization
            addressline={3-14-1, Hiyoshi, Kohoku-ku}, 
            city={Yokohama, Kanagawa},
            postcode={223-8522}, 
 %state={},
            country={Japan}}

%% Abstract
\begin{abstract}
It was recently demonstrated that the boundary element method based on the Burton--Miller formulation (BM-BEM), widely used for solving exterior problems, can be adapted to solve transmission problems efficiently. This approach utilises Calderon's identities to improve the spectral properties of the underlying integral operator. Consequently, most eigenvalues of the squared BEM coefficient matrix, i.e. the collocation-discretised version of the operator, cluster at a few points in the complex plane. When these clustering points are closely packed, the resulting linear system is well-conditioned and can be solved efficiently using the generalised minimal residual method with only a few iterations. However, when multiple materials with significantly different material constants are involved, some eigenvalues become separated, deteriorating the conditioning. To address this, we propose an enhanced Calderon-preconditioned BM-BEM with two strategies. First, we apply a preconditioning scheme inspired by the point Jacobi method. Second, we tune the Burton--Miller parameters to minimise the condition number of the coefficient matrix. Both strategies leverage a newly derived analytical expression for the eigenvalue clustering points of the relevant operator. Numerical experiments demonstrate that the proposed method, combining both strategies, is particularly effective for solving scattering problems involving composite penetrable materials with high contrast in material properties.
\end{abstract}


% \begin{highlights}
% \end{highlights}

%% Keywords
\begin{keyword}
Helmholtz' equation \sep
boundary element method (BEM) \sep
the Burton--Miller method \sep
transmission problem \sep
Calderon's preconditioning
\end{keyword}

\end{frontmatter}

\section{Introduction}
The boundary element method (BEM) is a numerical method that replaces the boundary value problem (BVP) of partial differential equations (PDEs) with boundary integral equations (BIEs). Unlike standard methods such as the finite element method (FEM) and finite-difference time-domain (FDTD) method, BEM solves BVP by discretising only the boundary of target objects. This feature makes BEM particularly suitable for wave scattering problems involving infinite domains, which are commonly encountered in acoustics, elasticity, and electromagnetics. Among such scattering problems, the so-called {\it transmission problem}, a mathematical model for wave scattering by penetrable objects, is of great engineering interest due to its wide range of potential applications~\cite{akduman_direct_2003, hammer_single_1998}.

So far, various boundary integral formulations for the transmission problem have been proposed, and some of them have been widely used in the BEM community. Among them, PMCHWT (Poggio--Miller--Chang--Harrington--Wu--Tsai) formulation~\cite{chew_fast_2001}, M\"uller formulation~\cite{muller_foundations_2013} and the single boundary integral equation (SBIE)~\cite{kleinman_single_1988, costabel1985direct} are especially well recognised as solvers that are free from real-valued fictitious eigenvalues. PMCHWT formulation, one of the most frequently used, can naturally be combined with {\it Calderon's preconditioning}~\cite{christiansen2000preconditionneurs, antoine_integral_2008}. This formulation can thus accelerate the convergence of iterative linear solvers for algebraic equations stemming from the BIE. Niino and Nishimura \cite{niino_preconditioning_2012} also found that the underlying (discretised version of) integral operator itself can have the role of (the inverse of) preconditioner by appropriately arranging the layer potentials, realising fast convergence for the generalised minimal residual method (GMRES)~\cite{saad_gmres_1986} without explicitly multiplying the preconditioner when the integral equation is discretised by the collocation. The method can easily be extended to the Galerkin BEM with only a few modifications with a cheap preconditioner. The M\"uller formulation provides an integral equation of the Fredoholm second type, thus providing a fast convergence for iterative solvers without any preconditioning. Both methods, however, may have complex-valued fictitious eigenvalues with a tiny imaginary part in its absolute, deteriorating the accuracy even at the real-valued incident frequency that is close to the complex-valued fictitious eigenvalue~\cite{misawa_boundary_2017}. On the other hand, the integral operator for SBIE does not involve fictitious eigenvalues around the real axis, thus achieving high efficiency and accuracy. The SBIE requires some special modifications in dealing with multiple materials, though~\cite{claeys2015second}. The multi-trace boundary integral formulation~\cite{hiptmair_multiple_2012} would also be an attractive choice, offering a straightforward extension to the multi-material problem. The formulation, however, requires multiple unknown densities on a single boundary element~\cite{hiptmair_multiple_2012}.

The Burton--Miller method~\cite{burton_application_1971} has the same distribution for the fictitious eigenvalues with the SBIE and is straightforwardly applicable to multi-material scattering. Iterative solvers for solving (the discretised version of) the boundary integral equations of Burton--Miller type (BM-BIE) can, however, be slow to convergent, especially in the case of transmission problems. Matsumoto et al \cite{matsumoto_calderon-preconditioned_2023} have recently proposed a Calderon preconditioning for the BM-BIE, in which they designed a well-conditioned integral operator in the sense that its square does not involve any hypersingular operator. Also, the eigenvalues of the square can only accumulate at several points in the complex plane. Especially for the transmission problems with a single material, there exists only one eigenvalue clustering point, reducing the iterative number for GMRES. This formulation can also be applied to multi-material problems. The key idea for accelerating multi-material cases is to use additional Burton--Miller equations not only for unbounded domains but also for some bounded domains. It has already been verified that the Calderon-preconditioned BM-BIE thus constructed achieves the faster convergence by GMRES than the conventional BM-BIE in many cases.  Nevertheless, when involving multiple materials with high contrasts in their material constants, such as permittivities and densities, the accumulation points can considerably separate with each other, and thus the GMRES convergence can become slow even with the preconditioning. 

Since it is crucial to provide an efficient solver for scattering problems involving composite structures with high-contrast material parameters, especially in complex situations where wave scattering is difficult to control, we present an improved version of the BM-BEM in this paper. We begin by revisiting the BIE formulation in \cite{matsumoto_calderon-preconditioned_2023}, showing that the eigenvalue accumulation can be analytically characterised in general settings. This insight allows us to construct a simple preconditioner, similar to the point Jacobi method, which collapses the clustered eigenvalues into a single point in the complex plane. We also identify both the optimal BM parameters and the interior subdomains to which the BM equations are applied, thereby achieving optimal conditioning of the resulting algebraic system. We demonstrate that the combination of these strategies yields the best overall performance. While the original study \cite{matsumoto_calderon-preconditioned_2023} demonstrated the efficiency of the approach only in two dimensions, we extend the analysis to three-dimensional problems, further confirming the versatility of both the original and improved Calderon-preconditioned BM-BEM.

The remainder of the paper is organised as follows. In Section 2, we review the BVP of interest and the associated BIEs combined with the Calderon preconditioning introduced in the previous study~\cite{matsumoto_calderon-preconditioned_2023}. We then present an analytical expression for the eigenvalue accumulation of the square of the relevant integral operator. Based on this result, Section 3 introduces an efficient point Jacobi-like preconditioning strategy, together with parameter optimisation techniques, to accelerate the BEM solver. Numerical examples in Section 4 demonstrate the performance of the proposed method, showing that it is highly efficient across a range of cases. Finally, Section 5 concludes the paper and discusses potential directions for future work.

\section{Formulation of BM-BIE with Calderon's preconditioning}\label{sec:BMBIE}
In this section, we examine in depth the recently proposed BEM incorporating Calderon's preconditioning~\cite{matsumoto_calderon-preconditioned_2023} as a preliminary step toward the presentation of the proposed point Jacobi-like preconditioning and parameter tuning. After stating the boundary value problem of interest in Section~\ref{subsec:bvp}, we outline the BIE formulated in the previous study in Section~\ref{subsec:matsumoto}. We then present an analytical investigation of the eigenvalue accumulation points of the underlying integral operator in Section~\ref{subsec:tomoyasu}. This theoretical insight will serve as a foundation for the developments in the subsequent sections.
\subsection{Statement of the boundary value problem}\label{subsec:bvp}
Let us consider that $\mathbb{R}^3$ is partitioned into $M$ subdomains as $\bigcup_{i=1}^M\overline{\Omega_i}$ with a sole unbounded domain $\Omega_1$. Given an incident field $u^{\mathrm{in}}:\Omega_1\rightarrow\mathbb{C}$ satisfying the three-dimensional Helmholtz equation of wavenumber $k_1$, the total field $u:\mathbb{R}^3\rightarrow\mathbb{C}$ solves the following boundary value problem:
\begin{align}
  \nabla^2u(\bs{x})+k_i^2u(\bs{x})=0\quad\bs{x}\in\Omega_i,
  \label{eq:helmholtz}\\
  \lim_{\eta\downarrow0}u(\bs{x}+\eta\bs{n})=\lim_{\eta\downarrow0}u(\bs{x}-\eta\bs{n})\quad\bs{x}\in\Gamma, 
  \label{eq:bcinu}\\
  \lim_{\eta\downarrow0}w(\bs{x}+\eta\bs{n})=\lim_{\eta\downarrow0}w(\bs{x}-\eta\bs{n})\quad\bs{x}\in\Gamma, 
  \label{eq:bcinw}\\
  \lim_{|\bs{x}|\rightarrow\infty}|\bs{x}|\left(\dfrac{\partial}{\partial|\bs{x}|}-\mathrm{i}k_1\right)\left(u(\bs{x})-u^{\mathrm{in}}(\bs{x})\right)=0,
  \label{eq:rc}
\end{align}
where $k_i>0$ denotes the wavenumber of waves propagating in $\Omega_i$, defined by $k_i:=\omega\sqrt{\varepsilon_i}$. Here, $\omega>0$ is the angular frequency of the incident wave, and $\varepsilon_i > 0$ (for $i = 1, \ldots, M$) is the material constant associated with $\Omega_i$. The boundary conditions given in (\ref{eq:bcinu}) and (\ref{eq:bcinw}) on the boundary $\Gamma:=\bigcup_{i=1}^M\partial\Omega_i$ impose jump relations on the total field and its flux. The flux here is defined as $w :=\frac{1}{\varepsilon_i} \bs{n} \cdot \nabla u$, where $\bs{n}$ is the unit normal vector on $\Gamma$. On $\partial \Omega_1$, $\bs{n}$ is oriented outward from $\Omega_1$. On the remaining parts of the boundary, the direction of $\bs{n}$ is chosen arbitrarily. We also define the boundary $\Gamma_{ij}:=\partial\Omega_i\cap\partial\Omega_j$, on which $\bs{n}$ is taken to point from $\Omega_i$ toward $\Omega_j$. Note that, under this convention, $\Gamma_{i1}$ (for $i = 2, \ldots, M$) does not exist. To ensure the uniqueness of the solution to the boundary value problem, the radiation condition (\ref{eq:rc}) is also imposed at infinity, where $\mathrm{i}$ denotes the imaginary unit.

\subsection{Calderon-preconditioned BM-BEM}\label{subsec:matsumoto}
To introduce the recently proposed Calderon-preconditioned BM-BIE, developed by Matsumoto et al~\cite{matsumoto_calderon-preconditioned_2023}, for a general setting, we first set up the necessary notation. Let ${\cal T}_i$ be the set of indices $j$ such that $\Omega_i$ and $\Omega_j$ share a common boundary, and let ${\cal T}_i^+$ (resp. ${\cal T}_i^-$) be a subset of ${\cal T}_i$, where the normal vector on $\partial \Omega_i \cap \partial \Omega_j$ is directed from $\Omega_i$ (resp. from $\Omega_j$ for each $j \in {\cal T}_i^-$). We also denote the traces of $u$ and $w$ on $\Gamma_{ij}$ by $u_{ij}$ and $w_{ij}$, respectively. With these notations, the boundary value problem \eqref{eq:helmholtz}--\eqref{eq:rc} is transformed into the following BIEs on $\bs{x}\in\Gamma_{ij}$: 
\begin{align}
-\frac{\alpha_1}{2}u_{ij}(\bs{x})
&
+\sum_{p\in{\cal T}_j^-}\alpha_1[{\cal D}^j_{\Gamma_{pj}} u_{pj}](\bs{x})
-\sum_{p\in{\cal T}_j^+}\alpha_1[{\cal D}^j_{\Gamma_{jp}} u_{jp}](\bs{x}) \notag \\
&
-\sum_{p\in{\cal T}_j^-}\alpha_1\varepsilon_j[{\cal S}^j_{\Gamma_{pj}} w_{pj}](\bs{x})
+\sum_{p\in{\cal T}_j^+}\alpha_1\varepsilon_j[{\cal S}^j_{\Gamma_{jp}} w_{jp}](\bs{x})=0, \label{eq:bie1}\\
\frac{1}{2}u_{ij}(\bs{x})
&
+\sum_{p\in{\cal T}_i^+}[({\cal D}^i_{\Gamma_{ip}}+\alpha_i{\cal N}^i_{\Gamma_{ip}})u_{ip}](\bs{x})
-\sum_{p\in{\cal T}_i^-}[({\cal D}^i_{\Gamma_{pi}}+\alpha_i{\cal N}^i_{\Gamma_{pi}})u_{pi}](\bs{x}) \notag\\
+\frac{\alpha_i\varepsilon_i}{2}&w_{ij}(\bs{x})
-\sum_{p\in{\cal T}_i^+}\varepsilon_i[({\cal S}^i_{\Gamma_{ip}}+\alpha_i({\cal D}^i_{\Gamma_{ip}})^*)w_{ip}](\bs{x})
+\sum_{p\in{\cal T}_i^-}\varepsilon_i[({\cal S}^i_{\Gamma_{pi}}+\alpha_i({\cal D}^i_{\Gamma_{pi}})^*)w_{pi}](\bs{x}) \notag \\
&=\delta_{i1}(u^\mathrm{in}_{1j}(\bs{x})+\alpha_1\varepsilon_1 w^\mathrm{in}_{1j}(\bs{x})), 
\label{eq:bie2}
\end{align}
 where $\delta_{ij}$ is the Kronecker delta, and $\alpha_i$ is the coupling coefficient of the Burton--Miller method. In this study, $\alpha_1=-\mathrm{i}/k_1$ is used to keep the complex-valued fictitious eigenvalue away from the real axis~\cite{zheng_is_2015}. The settings of $\alpha_i$ for $i\ne 1$ shall be discussed later. The integral operators in (\ref{eq:bie1}) and (\ref{eq:bie2}) are defined as
\begin{align}
  [\mathcal{S}^i_\Gamma w](\bs{x}):=\int_{\Gamma}G_i(\bs{x}-\bs{y})w(\bs{y})\mathrm{d}\Gamma(\bs{y}),
  \label{eq:slp}\\
  [\mathcal{D}^i_\Gamma u](\bs{x}):=\int_{\Gamma}\dfrac{\partial G_i(\bs{x}-\bs{y})}{\partial\bs{n}(\bs{y})}u(\bs{y})\mathrm{d}\Gamma(\bs{y}),
  \label{eq:dlp}\\
  [(\mathcal{D}^i_{\Gamma})^*w](\bs{x}):=\int_{\Gamma}\dfrac{\partial G_i(\bs{x}-\bs{y})}{\partial\bs{n}(\bs{x})}w(\bs{y})\mathrm{d}\Gamma(\bs{y}),
  \label{eq:nslp}\\
  [\mathcal{N}^i_\Gamma u](\bs{x}):=\mathrm{p.f.}\int_{\Gamma}\dfrac{\partial^2 G_i(\bs{x}-\bs{y})}{\partial\bs{n}(\bs{x})\partial\bs{n}(\bs{y})}u(\bs{y})\mathrm{d}\Gamma(\bs{y}),
  \label{eq:ndlp}
\end{align}  
where $G_i$ in (\ref{eq:slp})--(\ref{eq:ndlp}) is the following fundamental solution of the three-dimensional Helmholtz equation:
\begin{align}
  G_i(\bs{x}):=\dfrac{\mathrm{e}^{\mathrm{i}k_i|\bs{x}-\bs{y}|}}{4\pi|\bs{x}-\bs{y}|}
  \label{eq:green}
\end{align}
that satisfies the outgoing radiation condition, and ``p.f.'' in (\ref{eq:ndlp}) refers to the finite part of the diverging integral.

Let us now consider all index pairs $(i,j)$ such that $\Gamma_{ij} \ne \emptyset$. These pairs are collected into an ordered set, denoted by ${\cal B}$. The BIE system is then constructed by listing the equations corresponding to \eqref{eq:bie1} for all $(i,j) \in {\cal B}$ in order, followed by those corresponding to \eqref{eq:bie2} in the same order. Based on ${\cal B}$, we also define $u_{\cal B}$ and $w_{\cal B}$ as ordered lists of functions that collect all $u_{ij}$s and $w_{ij}$s, respectively, according to the order in ${\cal B}$. With these settings, we obtain the following system of BIEs:
\begin{align}
{\cal A}
\begin{pmatrix}
u_{\cal B}\\
w_{\cal B}
\end{pmatrix}
 = \text{right-hand side}, 
\label{eq:Ax=b}
\end{align}
where the integral operator $\cal A$ represents a block operator acting on the ordered pair $(u_{\mathcal{B}}, w_{\mathcal{B}})^t$. It is constructed by stacking the boundary integral operators corresponding to \eqref{eq:bie1} and \eqref{eq:bie2} for each $(i,j)\in \mathcal{B}$, in the order determined by $\mathcal{B}$. Accordingly, $\cal A$ consists of $2N_{\mathcal{B}} \times 2N_{\mathcal{B}}$ blocks, with $N_{\mathcal{B}}:=|\mathcal{B}|$. To clarify the notation and the procedure described above for constructing BIEs, we provide an illustrative example for a specific geometry in \ref{appendix:A}.

The system of BIEs \eqref{eq:Ax=b} differs from the standard formulation in two ways. First, the ``upper-half'' corresponding to \eqref{eq:bie1} is scaled by $-\alpha_1$. Second, all of the ``lower-half'' corresponding to \eqref{eq:bie2} are of the BM type. Just to avoid the fictitious eigenvalue problem, it is sufficient to apply the BM formulation only to the equation corresponding to the unbounded domain $\Omega_1$, i.e., the coupling coefficient $\alpha_i$ can be set to zero for all $i \ne 1$. Nonetheless, Matsumoto et al~\cite{matsumoto_calderon-preconditioned_2023} proposed setting nonzero values for all $\alpha_i$s red in (\ref{eq:bie2}). This, along with the first manipulation, improves the spectral properties of ${\cal A}$, as shall be discussed in the next subsection.

In the formulation of the BIE on $\Gamma_{ij}$, the trace from the domain $\Omega_j$ is handled by the standard integral equation \eqref{eq:bie1}, while that from $\Omega_i$ is handled by the BM equation \eqref{eq:bie2}. As mentioned above, however, in the case of $i \neq 1$, we may arbitrarily choose which trace to treat as the BM type. As shall be discussed later, the spectral properties of the integral operator $\mathcal{A}$ depend significantly on which side of the boundary is treated by the BM equation. In Section \ref{subsec:normal_flip}, we shall discuss how to appropriately select which integral equation to apply the Burton--Miller formulation to. Note that in our implementation, since we have decided to associate the orientation of the normal with the type of integral equation in a one-to-one correspondence as shown in \eqref{eq:bie1} and \eqref{eq:bie2}, the choice of which integral equation to apply the Burton--Miller formulation to can be made by adjusting the direction of the normal.

\subsection{Analytical evaluation of accumulated eigenvalues of the integral operator}\label{subsec:tomoyasu}
Let us now consider solving, by the GMRES~\cite{saad_gmres_1986}, the system of linear equations obtained by discretising the system of BIEs \eqref{eq:Ax=b} using the collocation method and accelerating the GMRES convergence. To this end, it is insightful to investigate the accumulation points of the eigenvalues of the square of the integral operator $\cal A$~\cite{niino_preconditioning_2012}. Through a detailed and rigorous calculation based on the integral equations \eqref{eq:bie1} and \eqref{eq:bie2}, and Calderon's formula, we found that ${\cal A}^2$ assumes the following form:
\begin{align}
  \mathcal{A}^2=\left(
  \begin{smallmatrix}
    \lambda_{1}\mathcal{I}&        &                                     &                                     &                 &             \\
                          & \ddots &                                     &                                     & \text{\huge{0}} &             \\
                          &        & \lambda_{{N_\mathcal{B}}}\mathcal{I}&                                     &                 &             \\
    c_{1}\mathcal{I}      &        &                                     & \lambda_{N_\mathcal{B}+1}\mathcal{I}&                 &             \\
                          & \ddots &                                     &                                     & \ddots          &             \\
    \text{\huge{0}}       &        & c_{N_\mathcal{B}}\mathcal{I}        &                                     &                 & \lambda_{2N_\mathcal{B}}\mathcal{I}
  \end{smallmatrix}  
  \right)+\mathcal{K}
  \label{eq:rep},
\end{align}
where $\cal I$ is the identity operator, $\lambda_i$ and $c_i$ are complex constants, and $\cal K$ represents an integral operator whose eigenvalues accumulate only at the origin of the complex plane~\cite{antoine_integral_2008,niino_preconditioning_2012}. The eigenvalue accumulation of the squared operator ${\cal A}^2$ is thus characterised by the operator whose upper-triangular blocks are all zero operators. The clustered eigenvalues $\lambda_{i_b}$ for $i_b=1, \ldots, 2N_{\mathcal{B}}$ are, of course, identical to the coefficients of the identity operators in the diagonal, given by
\begin{align}
  \left\{
  \begin{array}{ll}
    \lambda_{i_b}=\dfrac{\alpha_1^2}{4}\left(1+\dfrac{\alpha_i}{\alpha_1}\varepsilon_j\right)\\
    \lambda_{i_b+N_{\mathcal{B}}}=\dfrac{\alpha_1^2}{4}\left(\dfrac{\alpha_i^2}{\alpha_1^2}\varepsilon_i^2+\dfrac{\alpha_i}{\alpha_1}\varepsilon_j\right)
    \label{eq:param}
  \end{array}
  \right.,
\end{align}
for $1 \le i_b \le N_{\mathcal{B}}$, where $(i, j)$ is the $i_b$-th entry of the ordered list $\mathcal{B}$. The square of the integral operator in \eqref{eq:Ax=b} thus has at most $2N_{\mathcal{B}}$ accumulated eigenvalues, regardless of the choice of the arbitrary parameters $\alpha_i$ for $i=2,\ldots,M$. From \eqref{eq:param}, it is now clear why additional BM equations are necessary. When $\alpha_i = 0$, the operator acquires a zero eigenvalue, which can significantly increase the condition number of the BEM coefficient matrix obtained through discretisation.

Matsumoto et al~\cite{matsumoto_calderon-preconditioned_2023} set $\alpha_i=\alpha_1/\varepsilon_i$ based purely on empirical reasoning. By examining our findings in \eqref{eq:param}, however, one may observe that this setting actually leads to the following relation:
\begin{align}
 \lambda_{i_b}=\lambda_{i_b+N_\mathcal{B}}=\frac{\alpha_1^2}{4}\left(1+\frac{\varepsilon_j}{\varepsilon_i}\right)\quad\mathrm{for}\quad 1\le i_b\le N_\mathcal{B}.
\label{eq:lambda}
\end{align} 
This setting can therefore reduce the number of accumulation points to $N_\mathcal{B}$ from $2N_\mathcal{B}$. Most of the eigenvalues of the matrix $\mathsf{A}^2$ obtained by discretising the operator $\mathcal{A}^2$ via collocation are then expected to lie near these points, which may lead to faster convergence of GMRES. Note that it is not necessary to explicitly compute the squared matrix to take advantage of the improved conditioning of $\mathsf{A}^2$, since GMRES constructs approximate solutions in the Krylov subspace $\mathrm{span}\{\bs{x}_0,\mathsf{A}\bs{x}_0,\mathsf{A}^2\bs{x}_0, \ldots\}$, where $\bs{x}_0$ is the initial estimate for the solution of the system of linear equations.

The above $\alpha_i$ setting does, unfortunately, not always accelerate the GMRES convergence. For example, when the target domain consists of multiple materials with significantly different material constants, the eigenvalue accumulation points may be widely separated. In such cases, the condition number of $\mathsf{A}^2$ can still be large, resulting in poor convergence of GMRES. This motivates us to seek further improvements in the conditioning of the underlying system. In the next section, which presents the main contribution of this paper, we propose novel strategies to further improve the conditioning of the system.

\section{Strategies for further improving the BIE conditioning}\label{sec:main}
This section proposes two methods for the Calderon-preconditioned BM-BIE in the preceding section to improve its conditioning, particularly in cases where the contrast in material constants is large in multi-material settings. One method enforces eigenvalue accumulation at a single point by applying a simple preconditioner, while the other adjusts the parameters $\alpha_i$ for $i=2,\ldots, M$ and normal orientations to bring the eigenvalues closer together. 

\subsection{Point Jacobi-type preconditioning}\label{subsec:pj}
Since the accumulated eigenvalues \eqref{eq:param} for the squared integral operator ${\cal A}^2$ in the BIE system \eqref{eq:Ax=b} have now been analytically identified, one may use a point Jacobi-like preconditioner ${\cal M}$ to adjust the operator, such that all the clustered eigenvalues of $({\cal A}{\cal M}^{-1})^2$ are shifted to $1$. This can be achieved by defining $\cal M$ as
\begin{align}
{\cal M} = \text{diag} (\sqrt{\lambda_{1}}\mathcal{I}, \sqrt{\lambda_{2}}\mathcal{I}, \cdots, \sqrt{\lambda_{2N_{\mathcal{B}}}}\mathcal{I}),
\label{eq:precond}
\end{align}
where ${\cal M}$ has the same block structure as ${\cal A}$. It is thus expected that the collocation-discretised version of the (square of) preconditioned operator $\mathsf{({AM}^{-1})}^2$ is close to the identity matrix, with a considerably small condition number. As a result, we may expect an improvement in the convergence of the iterative solver.

It is worth noting that the preconditioning defined in this way can be both efficiently implemented and applied with little computational overhead. We can explicitly evaluate $\mathsf{M}^{-1}$ and store only the square root of the reciprocal of the estimated eigenvalues in \eqref{eq:param}, i.e., we only need an array of size equal to the number $2N_\mathcal{B}$ of the eigenvalue accumulation points, which hardly increases memory usage. In addition, we do not need to explicitly evaluate the preconditioned matrix $\mathsf{AM}^{-1}$. Instead, we simply multiply the diagonal matrix $\mathsf{M}^{-1}$ to the Krylov basis before multiplying by $\mathsf{A}$.

Note, however, that the squared coefficient matrix $\mathsf{A}^2$ may have eigenvalues that deviate from the estimated locations \eqref{eq:param} in the complex plane. If such an eigenvalue lies close to (resp. far from) the origin, multiplying the inverse of preconditioner $\mathsf{M}^{-1}$ with a tiny (resp. huge) diagonal entry may amplify the condition number of the preconditioned matrix $\mathsf{A}\mathsf{M}^{-1}$ and consequently that of its square $(\mathsf{A}\mathsf{M}^{-1})^2$. It is therefore important to construct $\mathsf{A}$ and $\mathsf{M}$ in such a way that $\mathsf{M}$ is close to the identity matrix, which corresponds to having all the eigenvalue accumulation points of $\mathcal{A}^2$ tightly packed. This may be achieved by appropriately setting the BIEs \eqref{eq:bie1} and \eqref{eq:bie2} with optimised tunable parameters $\alpha_i$ for $i = 2, \ldots, M$, which will be addressed in the next subsection.

\subsection{Minimising the maximum pairwise ratios of accumulated eigenvalues of ${\cal A}^2$}\label{seusec:psa}
In this subsection, we formulate the BIE system \eqref{eq:Ax=b} so as to minimise the maximum ratio between the absolute values of any two accumulated eigenvalues of the operator ${\cal A}^2$, while keeping the number of accumulated eigenvalues fixed at $N_\mathcal{B}$. Specifically, we consider which integral equations should be formulated in the BM form, and how the parameters $\alpha_i$ for $i=2,\ldots,M$ in the BM equations should be chosen. Recall that the former is achieved by appropriately setting the normal direction as discussed at the end of Subsection \ref{subsec:tomoyasu}. These approaches also improve the effectiveness of the point Jacobi method, as the eigenvalues become tightly clustered and $\mathcal{M}^{-1}$ in \eqref{eq:precond} becomes close to the identity operator.
\subsubsection{Possible $\alpha_i$ choices for $i=2,\ldots, M$}\label{subsec:alpha_choice}
Let us temporarily fix the normal directions and investigate all possible choices of the parameters $\alpha_i$ for $i=2,\ldots, M$ that lead to the eigenvalues of $\mathcal{A}^2$ to cluster at $N_\mathcal{B}$ points. The original setting 
\begin{align}
 \gamma_i:=\frac{\alpha_i}{\alpha_1}=\frac{1}{\varepsilon_i}
 \label{eq:pattern1}
\end{align}
in \cite{matsumoto_calderon-preconditioned_2023} satisfies this requirement. We shall henceforth refer to this setting as P1. Let us now take the $i_b$-th element $(i,j)$ and the $j_b$-th element $(k,\ell)$ from the index set $\mathcal{B}$, and consider the corresponding BM-BIEs \eqref{eq:bie2} defined on $\Gamma_{ij}$ and $\Gamma_{k\ell}$. We assume here that $\alpha_k$ is fixed according to P1, under which the relation $\lambda_{j_b} = \lambda_{j_b + N_\mathcal{B}}$ holds. With this setting, we may choose $\gamma_i$ such that $\lambda_{i_b} = \lambda_{j_b} = \lambda_{j_b + N_\mathcal{B}}$, yielding
\begin{align}
\gamma_i = \frac{\varepsilon_\ell}{\varepsilon_k \varepsilon_j}.
\label{eq:pattern2}
\end{align}
We refer to this setting as P2. It is also possible to set $\gamma_i$ such that $\lambda_{i_b+N_\mathcal{B}} = \lambda_{j_b} = \lambda_{j_b + N_\mathcal{B}}$, which gives the following P3: 
\begin{align}
\gamma_i = \frac{-\varepsilon_{j} + \sqrt{\varepsilon_{j}^2 + 4\varepsilon_i^2\left(1 + \varepsilon_{\ell} / \varepsilon_k\right)}}{2\varepsilon_i^2}.
\label{eq:pattern3}
\end{align}
In deriving \eqref{eq:pattern3}, we reject negative values of $\gamma_i$, since having $\alpha_i$ share the sign as $\alpha_1$ ensures that all fictitious and complex-valued eigenvalues arising from the BM equation have negative imaginary parts.

Recalling that the use of either P1 in \eqref{eq:pattern1}, P2 in \eqref{eq:pattern2}, or P3 in \eqref{eq:pattern3} for the artificial BM parameter effectively limits the number of accumulation points to $N_\mathcal{B}$, it is sufficient to choose the parameter that minimises the ratio between the absolute values of the accumulated eigenvalues of $\mathcal{A}^2$ among these options.

\begin{figure}[h]
  \begin{center}
   \includegraphics[width=4cm]{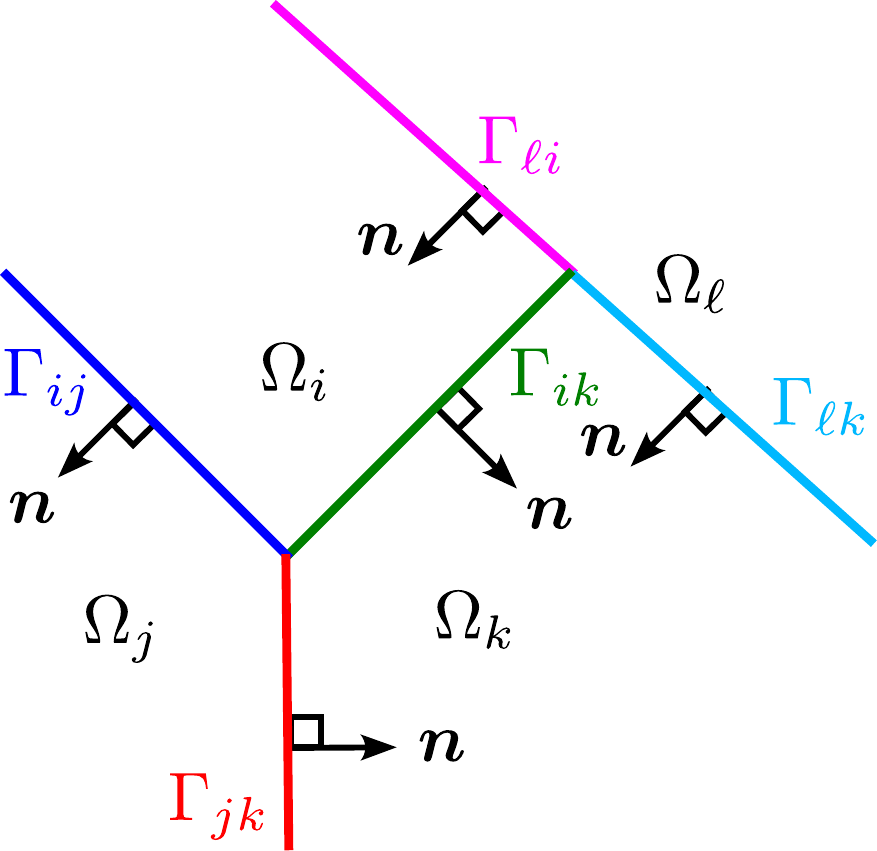}
   \caption{Illustration of condition C1, in which the parameter $\alpha_i$ is fixed using P1. Since the outward normal on $\partial \Omega_i$ points into both $\Omega_j$ and $\Omega_k$, the BM equations derived from $\Omega_i$ with $\alpha_i = \alpha_1 / \varepsilon_i$ are applied on $\Gamma_{ij}$ and $\Gamma_{ik}$. On $\Gamma_{\ell i}$, the BIE derived from $\Omega_i$ should be of the standard type, as the normal on this boundary points into $\Omega_i$. Instead, the BIE on $\Gamma_{\ell i}$ derived from $\Omega_\ell$ is of the BM type.}
   \label{fig:cond1}
  \end{center}
\end{figure}

In searching for the optimal parameter setting among the candidates (\ref{eq:pattern1})--(\ref{eq:pattern3}) while preserving the number of eigenvalue accumulation points, we must additionally impose the following two constraints:
\begin{itemize}
\item[\textbf{C1}] If multiple subdomains are adjacent to $\Omega_i$, and the outward normal vectors on $\partial \Omega_i$ point towards at least two of these adjacent subdomains, then $\gamma_i$ must be chosen according to P1 (\figref{fig:cond1}).
\item[\textbf{C2}] If $\varepsilon_i$ appears in $\gamma_k$ for the BM-BIE on $\Gamma_{k\ell}$ with $\ell \ne i$, then $\alpha_i$ must be fixed according to P1.
\end{itemize}

\subsubsection{Setting direction of normal vectors}\label{subsec:normal_flip}
\begin{figure}[h]
  \begin{center}
   \includegraphics[width=8cm]{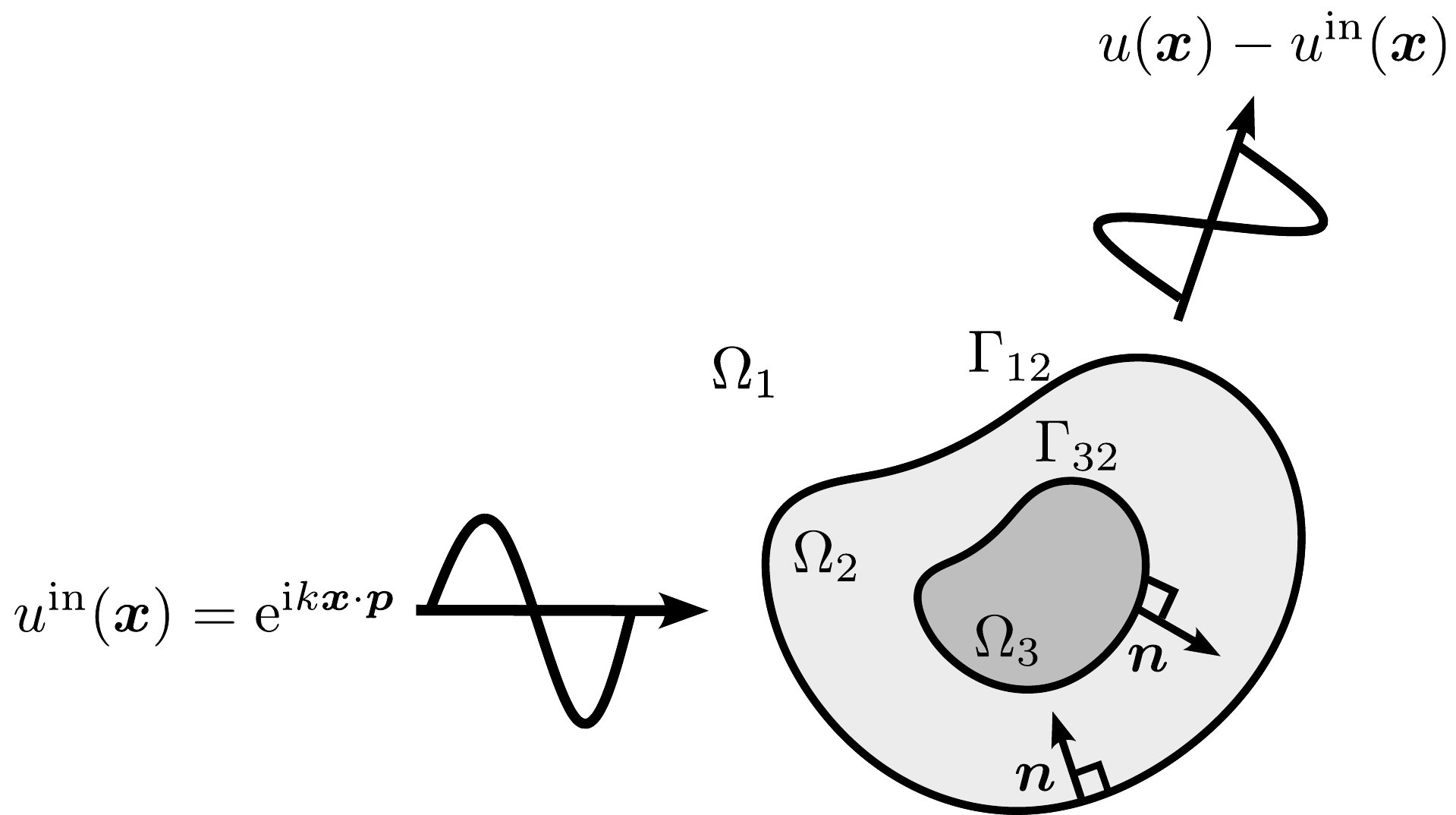}
   \caption{Multiple domains.}
   \label{fig:multi}
  \end{center}
\end{figure}

We then discuss an appropriate choice of normal directions, thereby ensuring that the suitable equations are formulated in the BM-type, as discussed in Subsection \ref{subsec:tomoyasu}, which improves the conditioning of the system of BIEs \eqref{eq:Ax=b}. To see how the normal direction setting influences the eigenvalue distribution of ${\cal A}^2$, let us consider the case where all the parameters $\alpha_i$s of P1 are used for a geometry depicted in Figure \ref{fig:multi}, with material constants $\varepsilon_1 = 1$, $\varepsilon_2 \simeq 1$, and $\varepsilon_3 \ll 1$.
In this case, $\mathcal{A}^2$ has two eigenvalue accumulation points, given by $\tfrac{\alpha_1^2(1+\varepsilon_2)}{4}$ and $\tfrac{\alpha_1^2(1+\varepsilon_2/\varepsilon_3)}{4}$, the ratio of whose absolute values is considerably large. On the other hand, when the normal vector on $\partial \Omega_2 \cap \partial \Omega_3$ is flipped, i.e. when $\Gamma_{32}$ is redefined as $\Gamma_{23}$, the accumulation points become $\tfrac{\alpha_1^2(1+\varepsilon_2)}{4}$ and $\tfrac{\alpha_1^2(1+\varepsilon_3/\varepsilon_2)}{4}$, which are much closer to each other than those in the previous configuration. Note that we cannot flip the normal on $\partial \Omega_1$, since the BIE corresponding to the unbounded domain $\Omega_1$ must be of the BM type to avoid the fictitious eigenvalue problem.

\subsubsection{Algorithm to optimise parameters $\alpha_i$ for $i=2,\ldots, M$}\label{subsec:algorithm}
By combining the parameter settings described in Subsection~\ref{subsec:alpha_choice} with normal direction flipping as discussed in Subsection~\ref{subsec:normal_flip}, we determine the optimal configuration of the BIE system \eqref{eq:Ax=b}. Specifically, we perform an exhaustive search over all admissible combinations of the parameter settings (\ref{eq:pattern1})--(\ref{eq:pattern3}), subject to C1 and C2, along with all possible normal orientations. For each configuration, we compute the accumulation points of the eigenvalues of ${\cal A}^2$ and evaluate the maximum ratio between the absolute values of any two such accumulation points. The configuration that minimises this maximum ratio is selected as the optimal one, as it yields the best spectral properties and thus improves the conditioning of the system.

\section{Numerical examples}
This section presents numerical experiments to verify that the proposed parameter tuning and simple preconditioning can accelerate the Calderon-preconditioned BM-BEM proposed  in~\cite{matsumoto_calderon-preconditioned_2023}. We first demonstrate in Subsection~\ref{sec:evaluate_matsumoto_method} that the conventional Calderon-preconditioned BM-BEM, without the proposed parameter tuning and point Jacobi-like preconditioning, performs well in many cases but may fail to accelerate GMRES convergence in others. We then show in Subsection~\ref{sec:effectiveness_of_the_proposed} that applying both strategies significantly improves the performance.

In all the following examples, we employ a collocation discretisation with piece-wise constant elements for the BM-BIEs. The singular parts of the boundary integrals are evaluated analytically, while the remaining parts are computed numerically using Gauss--Legendre (GL) quadrature with 16 nodes. Hypersingular integrals are regularised, and the resulting line integrals are evaluated using GL quadrature with 10 nodes. The resulting linear systems are solved using GMRES, with the tolerance set to $10^{-5}$.

\subsection{Performance evaluation of the conventional Calderon-preconditioned BM-BEM} \label{sec:evaluate_matsumoto_method}
\begin{figure}[h]
  \begin{center}
    \includegraphics[width=7cm]{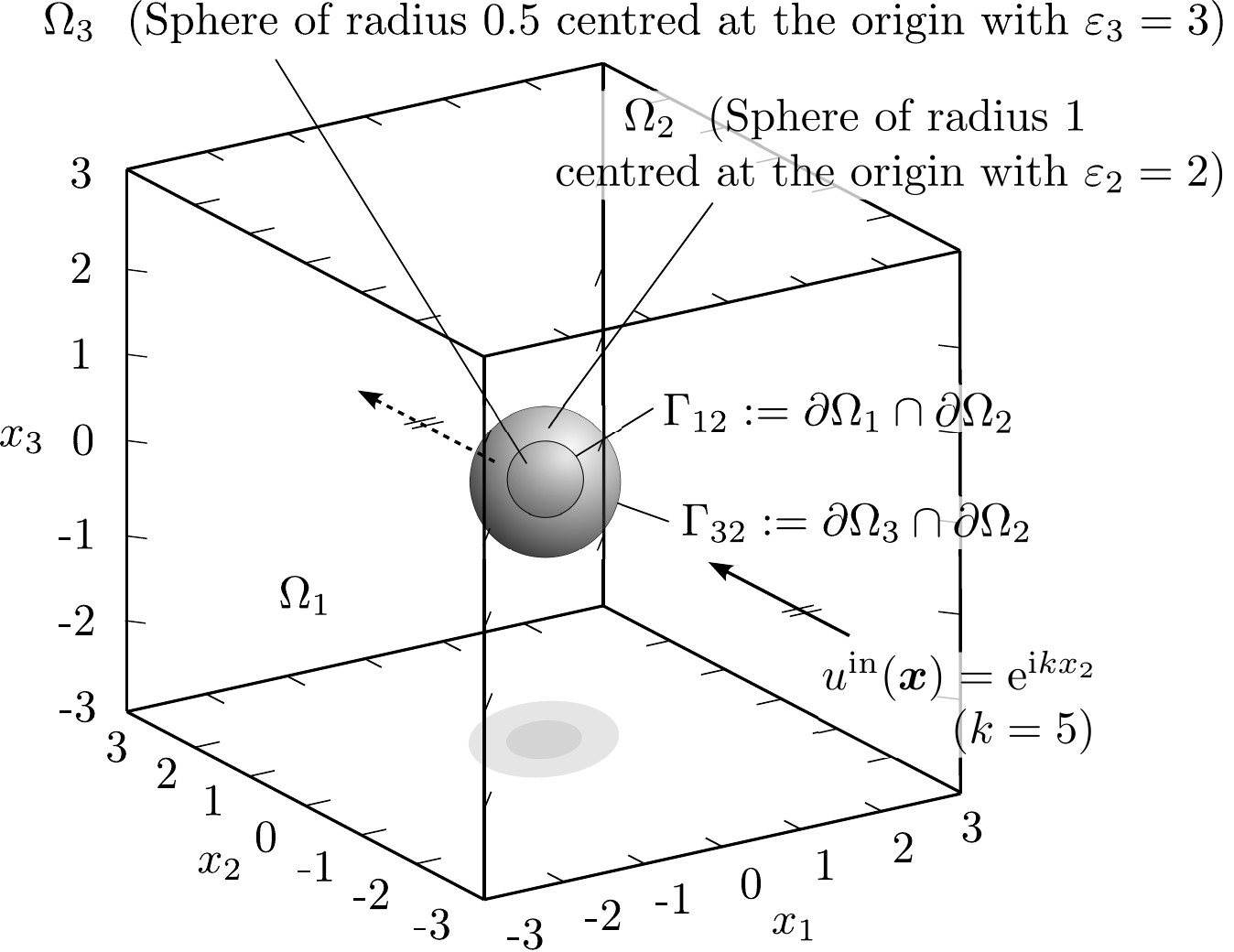}
    \caption{Problem setting for the transmission problem in the multi-material configuration. The normal vectors on the surfaces are directed from $\Omega_1$ and $\Omega_3$ in this setting. The boundaries are thus denoted as $\Gamma_{12}$ and $\Gamma_{32}$.}
    \label{fig:01kyu}
  \end{center}
\end{figure}

Let us first examine the accuracy and efficiency of the Calderon-preconditioned BM-BEM using a benchmark problem for which the analytical solution is available. To this end, we consider a multi-material configuration with $M = 3$, as illustrated in Figure \ref{fig:01kyu}. A spherical inclusion $\Omega_3$ of radius 0.5 is completely embedded within another spherical object $\Omega_2$ of radius 1.0, and both share the same centre at the origin. The outermost domain, the host matrix $\Omega_1$, along with the inclusions $\Omega_2$ and $\Omega_3$, have material constants 1, 2, and 3, respectively. As the incident wave, we use a plane wave propagating in the direction of the $x_2$-axis with angular frequency $\omega = 5.0$. 

In this configuration, according to our formulation, the unit normal vector on $\partial \Omega_1 \cap \partial \Omega_2$ should be consistently oriented outward from $\Omega_1$, since the BIE corresponding to the exterior domain $\Omega_1$ must be of BM type. On the other hand, the unit normal on $\partial \Omega_2 \cap \partial \Omega_3$ can, in principle, be arbitrarily oriented, but is here temporarily fixed to point outward from $\Omega_3$; that is, we apply the BM-BIE for $\Omega_3$ for the time being. As for the Burton--Miller coefficient $\alpha_3$, we temporarily fix it here according to P1, i.e., $\alpha_3 = \alpha_1 / \varepsilon_3$, as in~\eqref{eq:pattern1}. This choice is consistent with the original paper~\cite{matsumoto_calderon-preconditioned_2023}.

Here, we first check the relative $\ell_2$-error in $u$ at the collocation points on $\Gamma_{12} \cup \Gamma_{32}$, computed against the analytical solution, as well as the number of GMRES iterations required by the Calderon-preconditioned BM-BEM, with respect to the number of collocation points $N$. For comparison, we also compute the same quantities for the conventional BM-BEM, in which $\alpha_3 = 0$ and the integral operators in the upper \eqref{eq:bie1} and lower \eqref{eq:bie2} halves of the operator $\cal A$ in \eqref{eq:Ax=b} are interchanged. The $\ell_2$-error here is defined as
\begin{align}
  \mathrm{Error}:= 
    \sqrt{\frac{\sum_{i=1}^N |u_\mathrm{num}(\bs{x}_i)-u_\mathrm{ana}(\bs{x}_i)|^2}{\sum_{i=1}^N |u_\mathrm{ana}(\bs{x}_i)|^2}}
    \label{eq:error},
\end{align}
where $u_\mathrm{num}$ and $u_\mathrm{ana}$ are respectively the numerical and analytical solutions, $N$ is the number of the boundary elements, and $\bs{x}_i$ is the $i$-th collocation point, and the analytical solutions $u_\mathrm{ana}$ in (\ref{eq:error}) is derived from the following equations:
\begin{align}
  \left\{
  \begin{array}{ll}
  u_\mathrm{ana}(\bs{x})=u^\mathrm{in}(\bs{x})+\displaystyle\sum_{n=0}^\infty\displaystyle\sum_{m=-n}^n a_n^m h_n^{(1)}(k_1|\bs{x}|)Y_n^m(\theta,\phi)&\bs{x}\in\Omega_1\\
  u_\mathrm{ana}(\bs{x})=\displaystyle\sum_{n=0}^\infty\displaystyle\sum_{m=-n}^n(b_n^m h_n^{(1)}(k_2|\bs{x}|)+c_n^m j_n(k_2|\bs{x}|))Y_n^m(\theta,\phi)&\bs{x}\in\Omega_2\\
  u_\mathrm{ana}(\bs{x})=\displaystyle\sum_{n=0}^\infty\displaystyle\sum_{m=-n}^n d_n^m j_n(k_3|\bs{x}|)Y_n^m(\theta,\phi)&\bs{x}\in\Omega_3
  \label{eq:01ana3},
  \end{array}
  \right.
\end{align}
where $a_n^m,~b_n^m,~c_n^m$ and $d_n^m$ are complex constants determined, exploiting the orthogonality of the spherical harmonics, by the boundary conditions in (\ref{eq:bcinu}) and (\ref{eq:bcinw}). Also, $j_n$ and $h_n^{(1)}$ are respectively the $n^\mathrm{th}$ order spherical Bessel function and Hankel function of $1^\mathrm{st}$ kind, and $Y_n^m$ is spherical harmonics defined as
\begin{align}
  Y_n^m(\theta,\phi)=\sqrt{\frac{(n-m)!}{(n+m)!}}P_n^m(\cos\theta)\mathrm{e}^{\mathrm{i}m\phi}
  \label{eq:harmonic},
\end{align}
where $P_n^m$ is the associated Legendre polynomial, and $\theta$ and $\phi$ denote the polar and azimuthal angles for $\bs{x}$ respectively. In our computation, the infinite sums for $n$ in \eqref{eq:01ana3} is truncated by its first 51 terms. 

\begin{figure}[h]
  \begin{center}
    \includegraphics[scale=0.6]{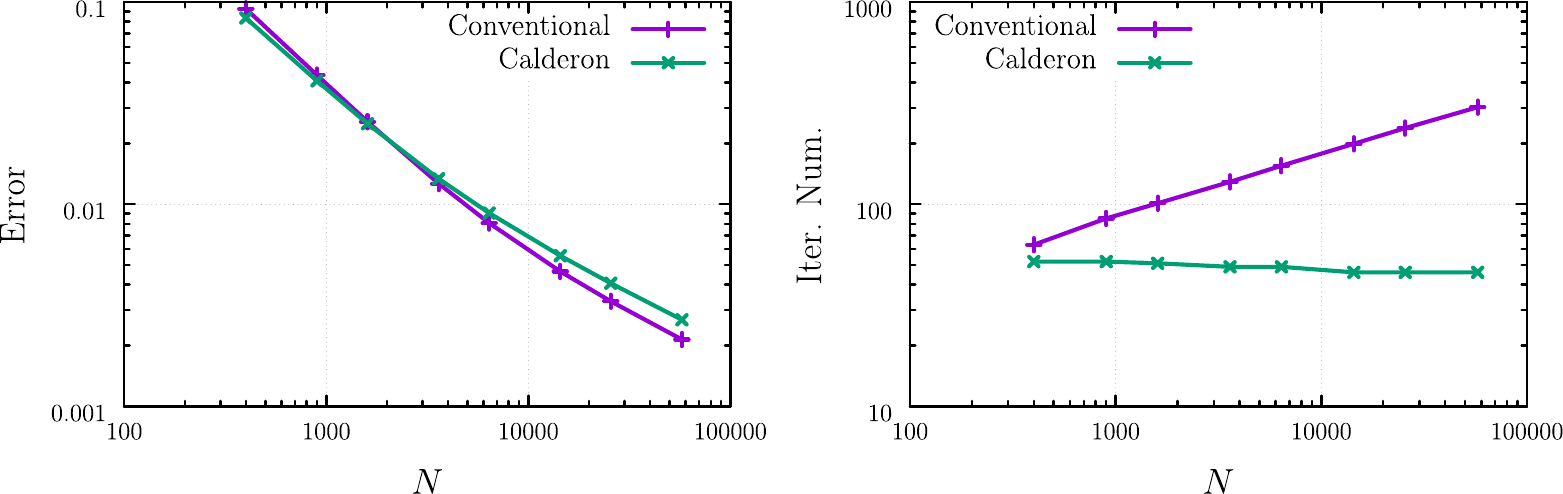}
    \caption{(Left) The relative $\ell_2$-error versus $N$ for wave scattering by concentric spheres. (Right) The number of GMRES iterations.}
    \label{fig:01elem-error-iter}
  \end{center}
\end{figure}

Figure~\ref{fig:01elem-error-iter} presents the results.
The left panel of Figure~\ref{fig:01elem-error-iter} shows that the accuracy of the preconditioned BM-BIE does not exactly match that of the conventional formulation, since the BM-BIE being solved differs between the two approaches. However, the difference in error is not significant. In the present case, the accuracy is slightly degraded for large $N$ by the proposed formulation, but this is not always the case; depending on the configuration, we have observed that the use of additional BM equations for interior domains can even improve accuracy in some cases~\cite{tomoyasu2024calderon}. The right panel shows that the Calderon preconditioning substantially reduces the number of GMRES iterations. More importantly, the proposed BM-BEM achieves a nearly constant iteration number regardless of the number of boundary elements $N$. \ref{appendix:B} provides further insight into how manipulating the eigenvalue distribution of the operator $\mathcal{A}^2$ leads to a reduced number of GMRES iterations.

\begin{figure}[h]
  \begin{center}
   \includegraphics[scale=0.63]{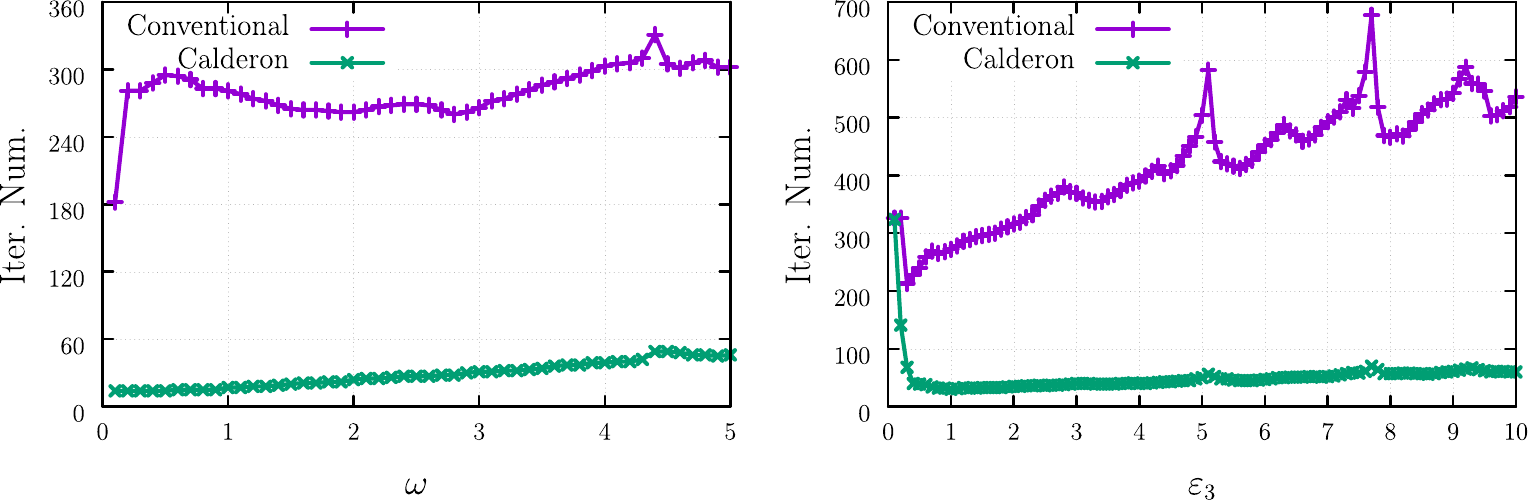}
   \caption{Number of GMRES iterations versus $\omega$ (left) and $\varepsilon_3$ (right) for wave scattering by concentric spheres.}
   \label{fig:01omega-epsln-iter}
  \end{center}
\end{figure}

To examine the versatility of the method, we evaluate the number of GMRES iterations under various settings of angular frequency $\omega$ and material constants $\varepsilon_i$ with the number of boundary elements fixed at $N = 57600$. The left panel of Figure~\ref{fig:01omega-epsln-iter} shows the case where the angular frequency $\omega$ is swept over the range $0 < \omega \le 5$, with the material constants fixed as $\varepsilon_i = i$ for $i = 1, \ldots, 3$. The right panel corresponds to the case where $\varepsilon_3$ is swept in the range $0 < \varepsilon_3 \le 10$, while the other parameters are set as $\varepsilon_1 = 1$ and $\varepsilon_2 = 1 / \varepsilon_3$. In this case, the angular frequency is fixed to $\omega = 5.0$.
Figure~\ref{fig:01omega-epsln-iter} demonstrates that the Calderon-preconditioned BEM requires fewer GMRES iterations than the conventional formulation for all tested parameter settings. It can, however, be observed that the performance may degrade significantly with a large contrast between $\varepsilon_2$ and $\varepsilon_3$. In particular, when $\varepsilon_3 \ll 1$ (and correspondingly $\varepsilon_2 \gg 1$), the number of iterations required by the preconditioned system is almost the same as that of the conventional one.

\begin{figure}[h]
  \begin{center}
   \includegraphics[width=7cm]{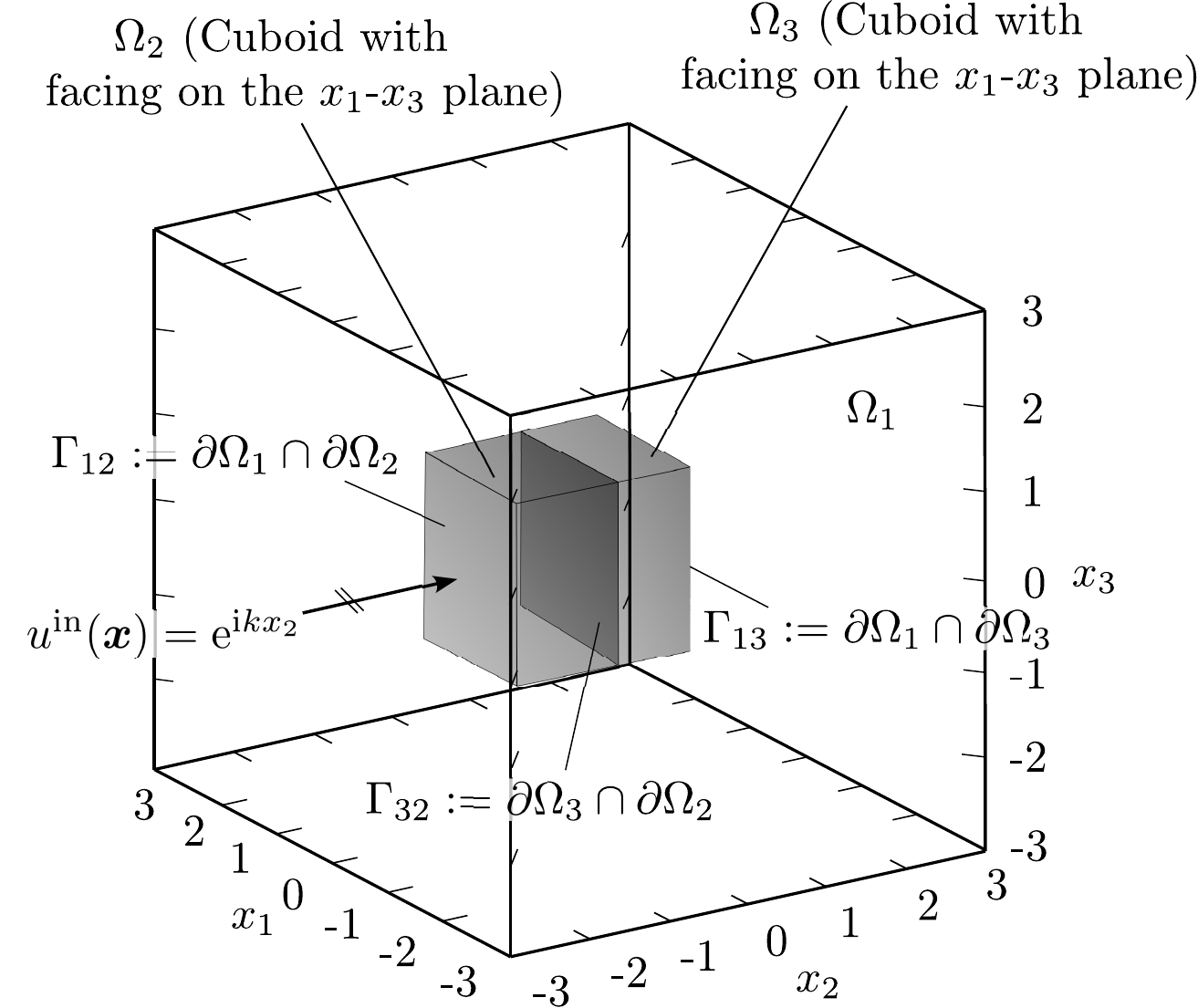}
   \caption{The problem setting for transmission problem with junctions and corners.}
   \label{fig:02cube}
  \end{center}
\end{figure}

We next examine a case involving geometry with junctions and corners. To this end, we consider the scattering by two cuboids, $\Omega_2$ and $\Omega_3$, which are in contact along one face, as illustrated in Figure~\ref{fig:02cube}.
\begin{figure}[h]
  \begin{center}
   \includegraphics[scale=0.42]{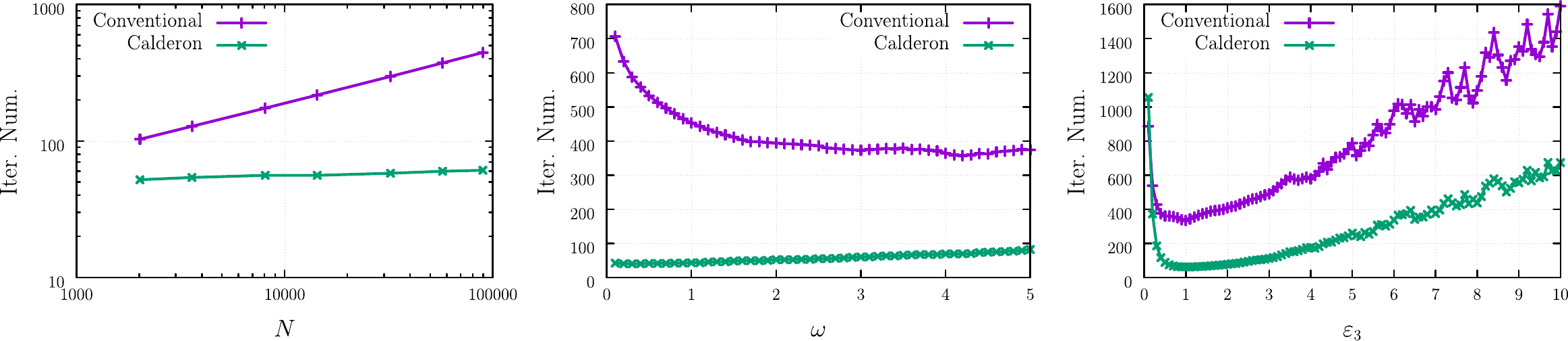}
   \caption{(Left) The number of GMRES iterations vs $N$ for wave scattering by a scatterer in Fig. \ref{fig:02cube} (centre) the number of iterations vs $\omega$ (right) the number of iterations vs $\varepsilon_2$.}
   \label{fig:02elem-omega-epsln-iter}
  \end{center}
\end{figure}
The left panel of Figure~\ref{fig:02elem-omega-epsln-iter} shows the number of GMRES iterations as a function of the number of boundary elements, $N$, using the same parameters as in the previous concentric spheres case; the angular frequency and material constants are fixed at $\omega = 5$, $\varepsilon_1=1$, $\varepsilon_2 = 2$, and $\varepsilon_3 = 3$. The centre and right panels correspond to those in Figure~\ref{fig:01omega-epsln-iter}, which show the dependence of the GMRES iterations on the angular frequency and material constants. For these computations, the number of boundary elements was fixed at $N = 57344$.

We can draw similar conclusions to those in the case of concentric spheres. The Calderon-preconditioned BM-BEM requires an almost constant number of GMRES iterations, regardless of the number of boundary elements, and generally performs better than the conventional method, even when the scatterers being analysed have junctions and corners. However, for material constants with a high contrast, both ends of the right panel in Figure~\ref{fig:02elem-omega-epsln-iter} show an increase in the number of GMRES iterations. In particular, when $\varepsilon_3 \ll 1$, the preconditioned BEM requires even more iterations than the conventional one. We therefore attempt to reduce the number of iterations in such high-contrast cases using the strategies described in Sections~\ref{sec:main}.

\subsection{Effectiveness of the proposed approaches} \label{sec:effectiveness_of_the_proposed}
In this subsection, we verify the effectiveness of the proposed point Jacobi-like preconditioning method and parameter tuning strategy. As an example, we revisit the concentric spheres case under the worst-case parameter settings shown in Figure~\ref{fig:01omega-epsln-iter}, namely, $\varepsilon_1 = 1.0$, $\varepsilon_2 = 0.1$, and $\varepsilon_3 = 10.0$. In the previous example shown in Figure~\ref{fig:01omega-epsln-iter}, the number of GMRES iterations was 324. We now apply the proposed point Jacobi-like preconditioner in Subsection~\ref{subsec:pj} while keeping the integral equation settings unchanged; that is, we adopt the BM-type BIE formulation for $\Omega_1$ and $\Omega_3$, and use the parameter $\alpha_3 = \alpha_1 / \varepsilon_3$ (referred to as P1) for the latter. The number of GMRES iterations is thereby reduced to 197. As previously noted, since the application of the inverse of the preconditioner matrix incurs negligible computational cost, this significant reduction represents a highly efficient improvement. 

While this alone is already promising, further gains can be achieved by tuning the formulation according to Subsection \ref{subsec:algorithm} prior to applying the preconditioner. In the current layered spherical geometry, we may choose the orientation of the interface $\partial \Omega_2 \cap \partial \Omega_3$, i.e. whether it is treated as $\Gamma_{32}$ or $\Gamma_{23}$, that determines which BIE is cast in BM form. Combined with the three parameter patterns (P1, P2, and P3), this gives rise to six possible BIE configurations.

The number of GMRES iterations for each of these configurations, with the point Jacobi-like preconditioning applied, is summarised in Table~\ref{table}.
\begin{table}[htbp]
  \centering
  \caption{Number of GMRES iterations for different combinations of interface orientation and parameter patterns, with the point Jacobi-like preconditioning applied. The value marked with ${}^*$ shows the untuned configuration discussed in the main text.}
  \label{table}
  \begin{tabular}{c|ccc}
    Interface & P1 & P2 & P3 \\
    \hline
    $\Gamma_{32}$ & 197${}^*$ & 181 & 181 \\
    $\Gamma_{23}$ & 219 & 296 & 151 \\
  \end{tabular}
\end{table}
The table highlights the importance of selecting an appropriate BIE system by properly orienting the interface and choosing suitable parameters before applying the preconditioner. 

\begin{figure}[h]
  \begin{center}
   \includegraphics[scale=0.6]{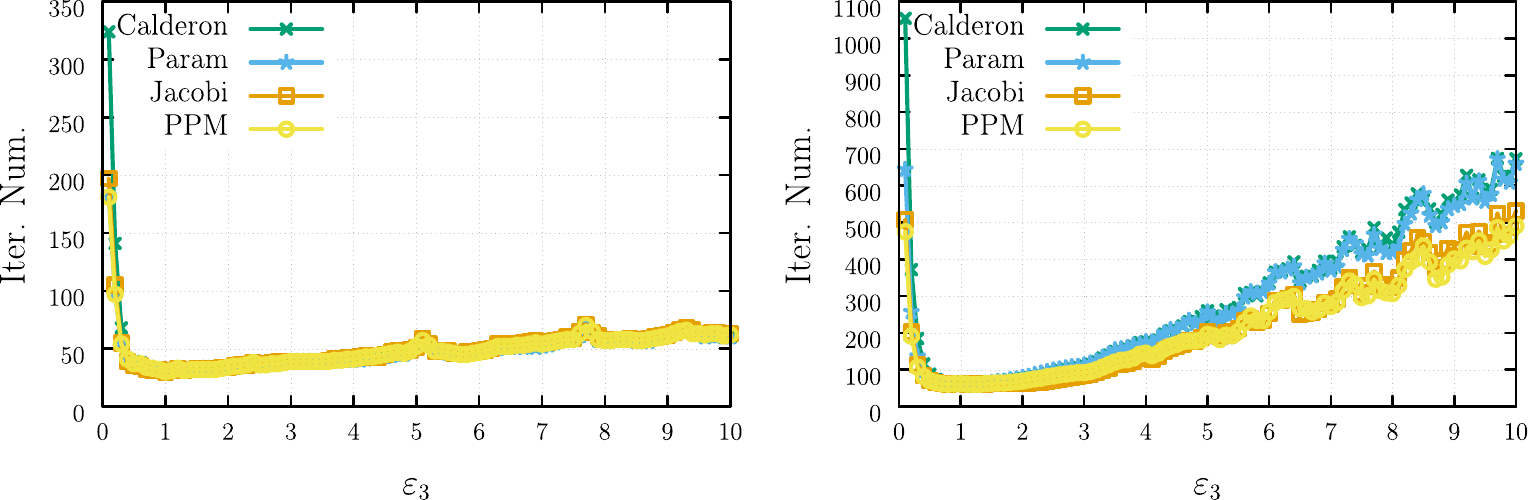}
   \caption{(Left) The number of GMRES iterations vs $\varepsilon_3$ in case of Fig. \ref{fig:01kyu} for wave scattering by an unit sphere (right) the number of iterations vs $\varepsilon_3$ in case of Fig. \ref{fig:02cube}.}
   \label{fig:0102epsln-iter}
  \end{center}
\end{figure}

We are thus motivated to combine the parameter tuning with point Jacobi-like preconditioning. Specifically, we first explore the best interface orientation and BM coefficient, aiming to minimise the maximum ratio of the absolute values of the accumulation points of the operator ${\cal A}^2$. We then solve the linear system using the point Jacobi-like preconditioner defined in \eqref{eq:precond}. With this strategy, we reproduce the results shown in the right panels of Figures~\ref{fig:01omega-epsln-iter} and~\ref{fig:02elem-omega-epsln-iter}, which indicate the number of GMRES iterations required for scattering problems involving concentric spheres and two cuboids, under various settings of the material constants. The results are summarised in Figure~\ref{fig:0102epsln-iter}. In the figure, the results obtained by using the proposed strategy (labelled as ``PPM''; Parameter-tuned and Point Jacobi-like Method) are compared those by other approaches. ``Calderon'' refers to the results presented in the previous subsection, obtained with a fixed interface orientation (as shown in Figures~\ref{fig:01kyu} and~\ref{fig:02cube}) and the BM parameter fixed to P1. ``Param'' indicates the use of parameter tuning alone, without any preconditioner, while ``Jacobi'' applies the point Jacobi preconditioner to the same linear system as ``Calderon''.
The figures demonstrate that the proposed method, which combines both strategies, improves efficiency in almost all cases, particularly in high-contrast material settings.

\begin{figure}[h]
  \begin{center}
   \includegraphics[width=8cm]{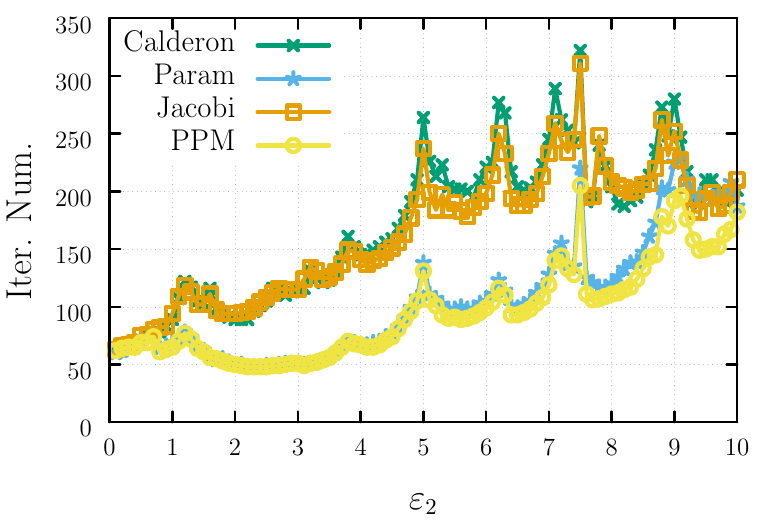}
   \caption{GMRES iteration number with varying material constant of $\Omega_2$, with $\varepsilon_1 = 1$ and $\varepsilon_3 = 10$ fixed.}
   \label{fig:01epsln2-iter}
  \end{center}
\end{figure}

Figure~\ref{fig:0102epsln-iter} may give the impression that the proposed method “PPM” performs almost identically to “Jacobi” in many cases, except in the extremely high-contrast settings. Further examination reveals, however, that this is not necessarily the case. To illustrate this, we again consider the scattering by concentric spheres shown in Figure~\ref{fig:01kyu}, with a different arrangement of material constants: $\varepsilon_1$ and $\varepsilon_3$ are fixed at $1$ and $10$, respectively, while $\varepsilon_2$ is varied from $0.1$ to $10$. The number of GMRES iterations required in this setting is presented in Figure~\ref{fig:01epsln2-iter}.
In the figure, over a wide range of $\varepsilon_2$, the performance of ``Jacobi'' is clearly observed to be inferior to that of ``PPM''. This result highlights the critical contribution of the proposed parameter tuning, which serves to enhance the effectiveness of ``Jacobi''. One also observes that, for $\varepsilon_2 > 9$, the performance of ``Param'' deteriorates significantly. These findings suggest that neither ``Jacobi'' nor ``Param'' alone is sufficient, and that the combined approach employed in ``PPM'' provides the most robust results.

\begin{figure}[h]
  \begin{center}
   \includegraphics[scale=0.6]{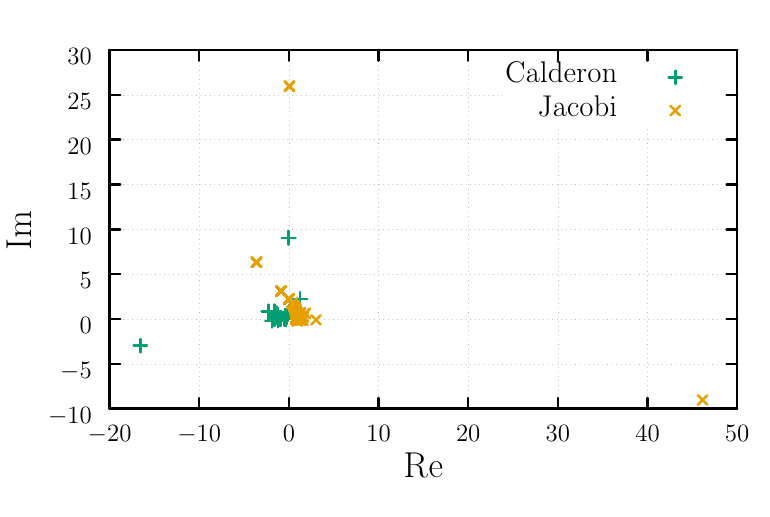}
   \includegraphics[scale=0.6]{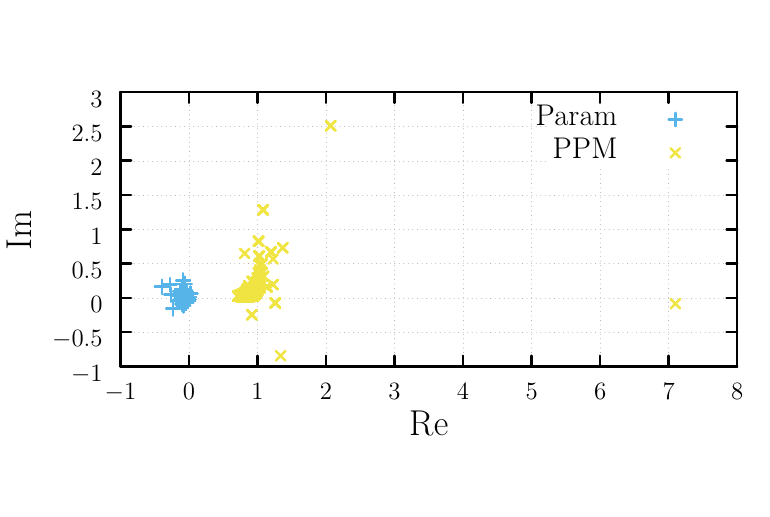}
   \caption{Change in the eigenvalue distribution of $\mathsf{A}^2$ from that of ``Calderon'' to  ``Jacobi'' (left) and ``Param'' to ``PPM'' (right). Note that the scales of the figures differ considerably.}
   \label{fig:01ev2}
  \end{center}
\end{figure}

Given that the eigenvalue accumulation points of the preconditioned squared integral operator $({\cal A}{\cal M}^{-1})^2$ are all theoretically forced to be at 1, the poor performance of ``Jacobi'' observed in Figure~\ref{fig:01epsln2-iter} may appear counter-intuitive. This discrepancy can be attributed to the fact that the collocation-discretised version of the operator, i.e., the square of preconditioned matrix $(\mathsf{A}\mathsf{M}^{-1})^2$ may deviate from the estimated accumulation point. To demonstrate this, we compute and compare the eigenvalue distributions of the squared matrices associated with both the ``Jacobi'' and ``PPM'' methods. The left panel of Figure~\ref{fig:01ev2} presents the spectrum for ``Jacobi'' along with that of ``Calderon'' method, while the right panel shows the result for ``PPM'' with ``Param''. In these computations, we set $\omega = 1$, $\varepsilon_1=1$, $\varepsilon_2 = 4$, $\varepsilon_3 = 10$, and used $N = 2880$ boundary elements.
The left panel of the figure shows that many of the eigenvalues for ``Calderon'' cluster around the expected locations, namely $-1.25$ and $-0.35$. These values are derived from equation~\eqref{eq:lambda} with $(i, j) = (1, 2)$ and $(3, 2)$, respectively. One can also observe, however, that a few of the eigenvalues lie significantly away from these points. In the eigenvalue distribution for “Jacobi”, also shown in the same panel, such outlying eigenvalues are further amplified in magnitude due to the structure of the preconditioner~\eqref{eq:precond}. Recall that the diagonal entries of  our preconditioner (its square to be precise) $(\mathsf{M}^{-1})^2$ correspond to the reciprocals of the cluster points. As a result, any deviation from the expected spectral location is magnified in this case. In contrast, this undesired effect can be mitigated in the proposed ``PPM'' method. Here, we first explore, via the normal flipping and parameter tuning, a modified BIE formulation whose eigenvalue accumulation points are as close together as possible before applying the preconditioning. The right panel of the figure shows the eigenvalue distribution for this formulation, labeled as ``Param''. The associated accumulation points are $-1.25$ and approximately $-1.0698$, the latter obtained from the first equation in~\eqref{eq:param} with $\alpha_i$ defined in~\eqref{eq:pattern3}, using the tuple $(i,j,k,l)=(2,3,1,2)$. With this design, the preconditioner also becomes closer to the identity matrix and does not scatter the outlier eigenvalues as strongly. Consequently, the eigenvalues in ``PPM'' tend to remain bounded in magnitude and well-behaved.

\begin{figure}[h]
\begin{center}
\includegraphics[scale=0.3]{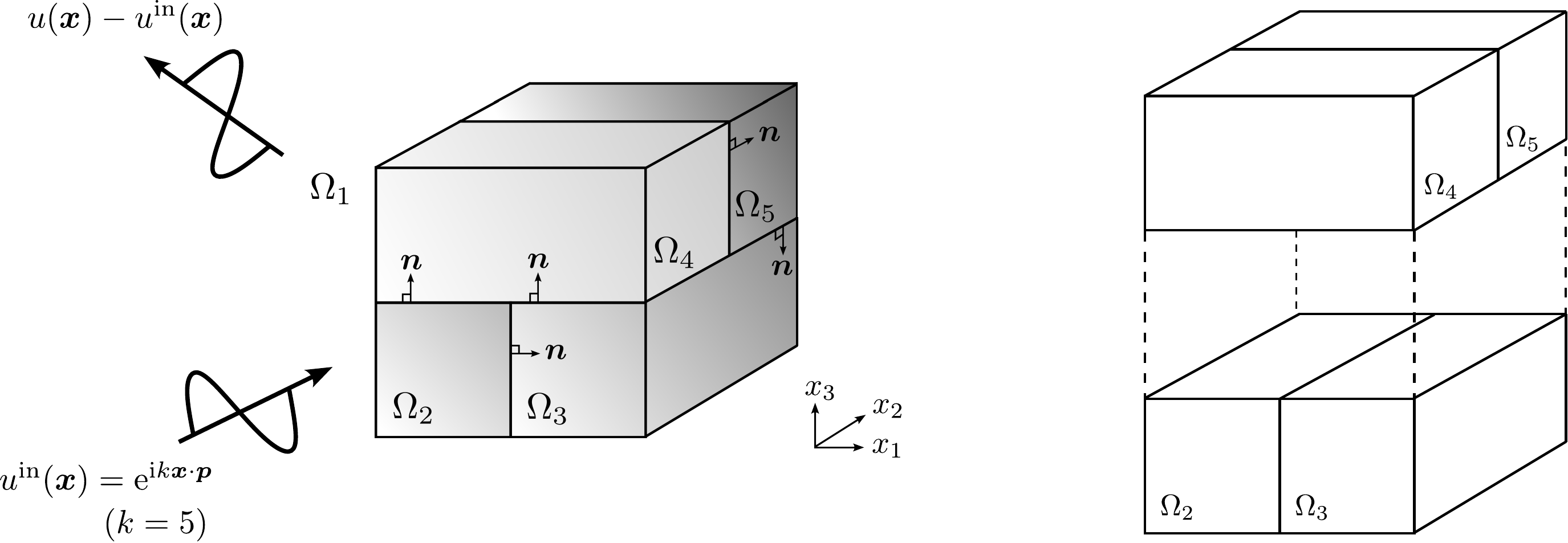}
\caption{Illustrative sketch of scattering by a composite scatterer consisting of four boxes (left), and its schematic diagram (right). The two lower boxes are aligned along the $x_2$-axis, while the two upper boxes are aligned along the $x_1$-axis.}
\label{fig:04cube}
\end{center}
\end{figure}

As the final numerical example, we examine the versatility of the proposed method in handling scattering from geometrically complex structures. To this end, we consider a composite scatterer composed of four stacked boxes, as shown in Figure~\ref{fig:04cube}.
Each box, as well as the host matrix, is characterised by specific material constants. Among these, some are fixed as $\varepsilon_1 = 1$, $\varepsilon_2 = 2$, and $\varepsilon_5 = 3$, while $\varepsilon_3$ and $\varepsilon_4$ vary: $\varepsilon_3$ ranges within $0 < \varepsilon_3 \le 10$, and $\varepsilon_4$ is set as $\varepsilon_4 = 1 / \varepsilon_3$. In the computation, the surfaces of the boxes are discretised into $N = 82944$ triangular boundary elements.

Since the structure is now considerably complex, many possibilities exist to construct the system of BIEs \eqref{eq:Ax=b}. To specify a particular system, we use the ordered set $\mathcal{B}$ of interface pairs, as defined in Subsection~\ref{subsec:matsumoto}. When the normal vectors on all interfaces are set as shown in Figure~\ref{fig:04cube}, this set is given by
\begin{align}
\mathcal{B} = \{ (1,2), (1,3), (1,4), (1,5), (2,3), (2,4), (3,4), (4,5), (5,2), (5,3) \}.
\end{align}
Following the formulation in Subsection~\ref{subsec:matsumoto}, we assign a BIE to each interface $\Gamma_{ij}$ with $(i,j)\in\mathcal{B}$ as follows:
\begin{itemize}
\item On the side of $\Omega_j$, the standard formulation in \eqref{eq:bie1} is used.
\item On the side of $\Omega_i$, the BM-type formulation in \eqref{eq:bie2} is employed.
\end{itemize}
The BM coefficient $\alpha_i$ used on $\Gamma_{ij}$ is thus set for the integral equation derived from $\Omega_i$. In particular, for those interfaces that do not involve the exterior domain $\Omega_1$, we allow multiple choices for $\alpha_i$, corresponding to P1, P2, or P3. To simplify the notation, we introduce $\Gamma_{ij}(\mathrm{P}^k)$ to indicate that:
\begin{itemize}
\item the interface is oriented from $\Omega_i$ to $\Omega_j$,
\item the BM-type BIE is applied to $\Omega_i$ with parameter choice P$k$ (e.g., $\alpha_i = \alpha_1 / \varepsilon_i$ for $k=1$),
\item and the standard BIE is applied to $\Omega_j$.
\end{itemize}
With this convention, the entire system of BIEs is fully characterised by the set of such interface notations. For instance, if all interfaces are oriented as shown in the figure, and all BM coefficients are set using P1, the BIE system can be denoted as:
\begin{align}
\Gamma_{23}(\mathrm{P}^1),\quad \Gamma_{24}(\mathrm{P}^1),\quad \Gamma_{34}(\mathrm{P}^1),\quad \Gamma_{45}(\mathrm{P}^1),\quad \Gamma_{52}(\mathrm{P}^1),\quad \Gamma_{53}(\mathrm{P}^1).
\label{eq:04cube_calderon}
\end{align}
Here, the interfaces $\Gamma_{1j}$ for all $j$ are omitted, as the corresponding BIEs have already been fixed as the preceding discussion. Note that for the parameter $\alpha_i$ on $\Gamma_{ij}$ with settings P2 or P3, we need to specify an additional element $(k, \ell)\in\mathcal{B}$, as seen in \eqref{eq:pattern2} and \eqref{eq:pattern3}. We shall thus use the notation $\Gamma_{ij}(\mathrm{P}^{2}_{k\ell})$ for this purpose.

\begin{figure}[h]
\begin{center}
\includegraphics[scale=0.7]{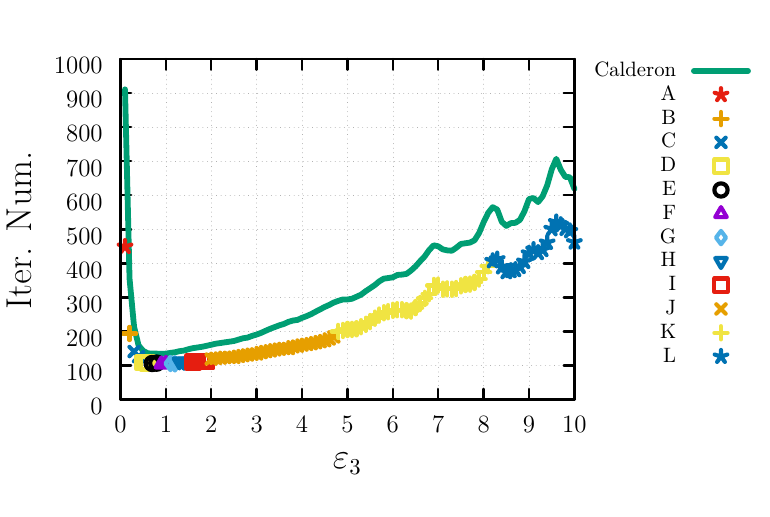}
\caption{The number of GMRES iterations corresponding to the case in Figure~\ref{fig:04cube}.}
\label{fig:04epsln-iter}
\end{center}
\end{figure}

Figure~\ref{fig:04epsln-iter} shows the number of GMRES iterations versus $\varepsilon_3$, obtained by the proposed method that combines boundary orientation flipping, parameter tuning, and point Jacobi-type preconditioning.
The labels A to L in the figure indicate the specific settings adaptively selected by the proposed method. Each of these configurations is detailed in Table \ref{table2}.
\begin{table}
\centering
\caption{List of interface orientations and BM-type parameter patterns selected by the proposed strategy for each configuration A–L in the scattering problem illustrated in Figure~\ref{fig:04cube}.}\label{table2}
 \begin{tabular}{l| c c c c c c}
         &  &  &  &  &  &  \\
    \cline{1-7}
    A  & $\Gamma_{23}(\mathrm{P}^1)$ & $\Gamma_{24}(\mathrm{P}^1)$ & $\Gamma_{34}(\mathrm{P}^2_{12})$ & $\Gamma_{45}(\mathrm{P}^2_{52})$ & $\Gamma_{52}(\mathrm{P}^1)$ & $\Gamma_{53}(\mathrm{P}^1)$ \\
    B  & $\Gamma_{23}(\mathrm{P}^1)$ & $\Gamma_{24}(\mathrm{P}^1)$ & $\Gamma_{43}(\mathrm{P}^3_{14})$ & $\Gamma_{54}(\mathrm{P}^1)$ & $\Gamma_{52}(\mathrm{P}^1)$ & $\Gamma_{35}(\mathrm{P}^2_{12})$ \\
    C  & $\Gamma_{23}(\mathrm{P}^1)$ & $\Gamma_{24}(\mathrm{P}^1)$ & $\Gamma_{43}(\mathrm{P}^2_{23})$ & $\Gamma_{54}(\mathrm{P}^1)$ & $\Gamma_{52}(\mathrm{P}^1)$ & $\Gamma_{35}(\mathrm{P}^2_{12})$ \\
    D  & $\Gamma_{32}(\mathrm{P}^2_{12})$ & $\Gamma_{24}(\mathrm{P}^1)$ & $\Gamma_{43}(\mathrm{P}^1)$ & $\Gamma_{45}(\mathrm{P}^1)$ & $\Gamma_{25}(\mathrm{P}^1)$ & $\Gamma_{53}(\mathrm{P}^2_{43})$ \\
    E  & $\Gamma_{32}(\mathrm{P}^2_{12})$ & $\Gamma_{24}(\mathrm{P}^1)$ & $\Gamma_{43}(\mathrm{P}^1)$ & $\Gamma_{45}(\mathrm{P}^1)$ & $\Gamma_{25}(\mathrm{P}^1)$ & $\Gamma_{53}(\mathrm{P}^3_{15})$ \\
    F  & $\Gamma_{32}(\mathrm{P}^2_{12})$ & $\Gamma_{42}(\mathrm{P}^1)$ & $\Gamma_{43}(\mathrm{P}^1)$ & $\Gamma_{45}(\mathrm{P}^1)$ & $\Gamma_{25}(\mathrm{P}^2_{12})$ & $\Gamma_{53}(\mathrm{P}^3_{15})$ \\
    G  & $\Gamma_{32}(\mathrm{P}^2_{12})$ & $\Gamma_{42}(\mathrm{P}^1)$ & $\Gamma_{43}(\mathrm{P}^1)$ & $\Gamma_{45}(\mathrm{P}^1)$ & $\Gamma_{25}(\mathrm{P}^2_{12})$ & $\Gamma_{53}(\mathrm{P}^3_{45})$ \\
    H  & $\Gamma_{23}(\mathrm{P}^1)$ & $\Gamma_{42}(\mathrm{P}^2_{12})$ & $\Gamma_{34}(\mathrm{P}^1)$ & $\Gamma_{54}(\mathrm{P}^3_{15})$ & $\Gamma_{25}(\mathrm{P}^1)$ & $\Gamma_{35}(\mathrm{P}^1)$ \\
    I  & $\Gamma_{23}(\mathrm{P}^1)$ & $\Gamma_{42}(\mathrm{P}^2_{12})$ & $\Gamma_{34}(\mathrm{P}^1)$ & $\Gamma_{54}(\mathrm{P}^2_{34})$ & $\Gamma_{25}(\mathrm{P}^1)$ & $\Gamma_{35}(\mathrm{P}^1)$ \\
    J  & $\Gamma_{23}(\mathrm{P}^1)$ & $\Gamma_{24}(\mathrm{P}^1)$ & $\Gamma_{34}(\mathrm{P}^2_{24})$ & $\Gamma_{45}(\mathrm{P}^2_{12})$ & $\Gamma_{52}(\mathrm{P}^1)$ & $\Gamma_{53}(\mathrm{P}^1)$ \\
    K  & $\Gamma_{23}(\mathrm{P}^1)$ & $\Gamma_{24}(\mathrm{P}^1)$ & $\Gamma_{43}(\mathrm{P}^2_{12})$ & $\Gamma_{54}(\mathrm{P}^2_{24})$ & $\Gamma_{25}(\mathrm{P}^1)$ & $\Gamma_{35}(\mathrm{P}^3_{12})$ \\
    L  & $\Gamma_{23}(\mathrm{P}^1)$ & $\Gamma_{24}(\mathrm{P}^1)$ & $\Gamma_{43}(\mathrm{P}^2_{12})$ & $\Gamma_{54}(\mathrm{P}^2_{14})$ & $\Gamma_{25}(\mathrm{P}^1)$ & $\Gamma_{35}(\mathrm{P}^3_{12})$ \\
\end{tabular}
\end{table}
As can be seen, the interface orientations and parameter patterns vary significantly across the different configurations, indicating the necessity of adaptively selecting both in order to optimise convergence. For reference, we also plot the result obtained using the setting fixed as in \eqref{eq:04cube_calderon}, labelled as Calderon. It is immediately evident that the proposed approach achieves superior performance for all tested cases.

\section{Conclusion}
This study proposed an extension to the recently developed Calderon-preconditioned BM-BEM~\cite{matsumoto_calderon-preconditioned_2023}. Although the original formulation ensures that the square of the underlying operator has only a few eigenvalue accumulation points with finite absolute values, these eigenvalues may still be widely spread, thereby deteriorating the conditioning of the resulting algebraic system in the BEM.

To address this issue, we first tuned the relevant parameters to cluster the accumulation points as tightly as possible. We then applied a point Jacobi-like preconditioning to further enhance the clustering. Both of these steps are made possible by having analytically identified the eigenvalue accumulation points of the relevant operator in a general setting. The resulting formulation achieves improved efficiency in the BM-BIE for transmission problems involving multiple materials with significantly different material constants.

We plan to extend the present approach to transmission problems in elastodynamics and electromagnetics in 3D. It would also be of interest to apply the proposed BM-BEM to frequency-differentiated Helmholtz' equation for wide-band frequency sweeps~\cite{chen2023generalized, kook2013acoustical}. Furthermore, we intend to explore a range of applications, including topology optimisation~\cite{nakamoto_fast_2016, isakari_topology_2014}.

\appendix
\setcounter{figure}{0}
\section{Calderon preconditioned BM-BIE for a specific geometry}\label{appendix:A}
\begin{figure}[h]
  \begin{center}
   \includegraphics[scale=0.22]{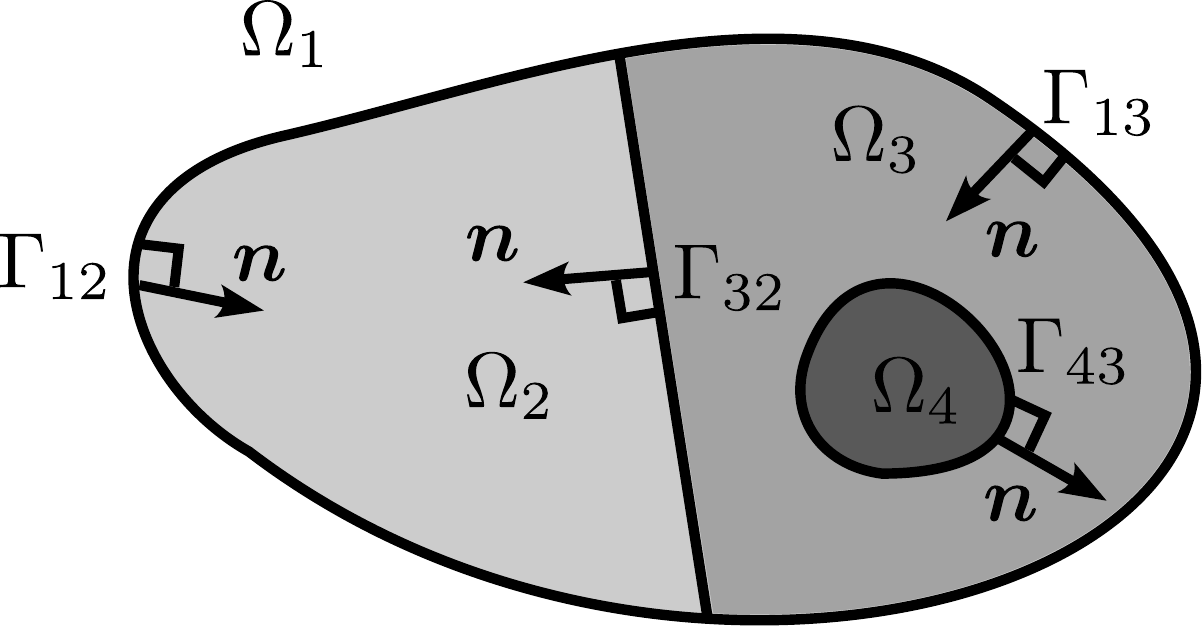}
   \caption{Composite material with 4 subdomains. The normal direction on $\Gamma_{ij} := \partial \Omega_i \cap \partial \Omega_j$ is defined to point outward from $\Omega_i$.}
   \label{fig:4domains}
  \end{center}
\end{figure}

As an example of the formulation presented in Section \ref{subsec:matsumoto}, we here derive the Calderon-preconditioned BIE system corresponding to wave scattering by the composite material shown in Figure~\ref{fig:4domains}.

Let us first consider, for example, the index set ${\cal T}_3$. We have ${\cal T}_3 = \{1, 2, 4\}$, since $\Omega_3$ shares boundaries with $\Omega_1$, $\Omega_2$, and $\Omega_4$. The subsets ${\cal T}_3^+ = \{2\}$ and ${\cal T}_3^- = \{1, 4\}$ follow from the orientation of the normal vectors: the normal on $\Gamma_{32}$ is defined to point outward from $\Omega_3$, whereas those on $\Gamma_{13}$ and $\Gamma_{43}$ point inward to $\Omega_3$. We then collect all index pairs $(i, j)$ such that $\Gamma_{ij}$ exists into an ordered set, denoted by $\mathcal{B} = \{(1,2), (1,3), (3,2), (4,3)\}$. The ordering of the pairs can be chosen arbitrarily, but must remain fixed once specified. In accordance with the order specified by $\mathcal{B}$, the traces $u_{ij}$ and $w_{ij}$ are collected into the sequences $u_{\mathcal{B}} = (u_{12}, u_{13}, u_{32}, u_{43})^t$ and $w_{\mathcal{B}} = (w_{12}, w_{13}, w_{32}, w_{43})^t$, respectively. The BIEs \eqref{eq:bie1}, derived as the trace of the integral representation in $\Omega_j$ to $\Gamma_{ij}$, are then listed in the order specified by $\mathcal{B}$ as follows:
\begin{align}
-\alpha_1\left(\frac{\cal I}{2}-{\cal D}^2_{12}\right)u_{12}
+\alpha_1{\cal D}^2_{32}u_{32}
-\alpha_1\varepsilon_2{\cal S}^2_{12}w_{12}
-\alpha_1\varepsilon_2{\cal S}^2_{32}u_{32}&=0,\label{eq:A1} \\
-\alpha_1\left(\frac{\cal I}{2}-{\cal D}^3_{13}\right)u_{13}
-\alpha_1{\cal D}_{32}^3u_{32}
+\alpha_1{\cal D}_{43}^3u_{43}
-\alpha_1\varepsilon_3{\cal S}^3_{13}w_{13}
+\alpha_1\varepsilon_3{\cal S}^3_{32}w_{32}
-\alpha_1\varepsilon_3{\cal S}^3_{43}w_{43}
&=0,\\
\alpha_1{\cal D}^2_{12}u_{12}
-\alpha_1\left(\frac{\cal I}{2}-{\cal D}^2_{32}\right)u_{32}
-\alpha_1\varepsilon_2{\cal S}^2_{12}w_{12}
-\alpha_1\varepsilon_2{\cal S}^2_{32}w_{32}
&=0,\\
\alpha_1{\cal D}^3_{13}u_{13}
-\alpha_1{\cal D}^3_{32}u_{32}
-\alpha_1\left(\frac{\cal I}{2}-{\cal D}^3_{43}\right)u_{43}%\notag \\
-\alpha_1\varepsilon_3{\cal S}^3_{13}w_{13}
+\alpha_1\varepsilon_3{\cal S}^3_{32}w_{32}
-\alpha_1\varepsilon_3{\cal S}^3_{43}w_{43}
&=0, 
\end{align}
in which all the equations are multiplied by the same constant $-\alpha_1$.
Subsequently, the BIEs \eqref{eq:bie2} are  listed, again in accordance with the order $\mathcal{B}$ as
\begin{align}
\left(\frac{\cal I}{2}+{\cal W}^1_{12}\right)u_{12}
+{\cal W}^1_{13}u_{13}
+\varepsilon_1\left(\frac{\alpha_1{\cal I}}{2}-{\cal V}^1_{12}\right)w_{12}
-\varepsilon_1{\cal V}^1_{13}w_{13}
=0,\\
{\cal W}^1_{12}u_{12}
+\left(\frac{\cal I}{2}+{\cal W}^1_{13}\right)u_{13}
-\varepsilon_1{\cal V}^1_{12}w_{12}
+\varepsilon_1\left(\frac{\alpha_1{\cal I}}{2}-{\cal V}^1_{13}\right)w_{13}
=0,\\
-{\cal W}^3_{13}u_{13}
+\left(\frac{\cal I}{2}+{\cal W}^3_{32}\right)u_{32}
-{\cal W}^3_{43}u_{43}
+\varepsilon_3{\cal V}^3_{13}w_{13}
+\varepsilon_3\left(\frac{\alpha_3{\cal I}}{2}-{\cal V}^3_{32}\right)w_{32}
+\varepsilon_3{\cal V}^3_{43}w_{43}
=0,\\
\left(\frac{\cal I}{2}+{\cal W}^4_{43}\right)u_{43}
+\varepsilon_4\left(\frac{\alpha_4{\cal I}}{2}-{\cal V}^4_{43}\right)w_{43}
=0,\label{eq:A8}
\end{align}
where the BM operators ${\cal V}^i_{jk}:={\cal S}^i_{jk}+\alpha_i({\cal D}^i_{jk})^*$ and ${\cal W}^i_{jk}:={\cal D}^i_{jk}+\alpha_i{\cal N}^i_{jk}$ are defined. Then, all the BIEs \eqref{eq:A1}--\eqref{eq:A8} are assembled in the specified order to form the whole system \eqref{eq:Ax=b}.

\setcounter{figure}{0}
\section{Eigenvalue distribution of $\mathcal{A}^2$}\label{appendix:B}

\begin{figure}[h]
   \begin{center}
    \includegraphics[scale=0.63]{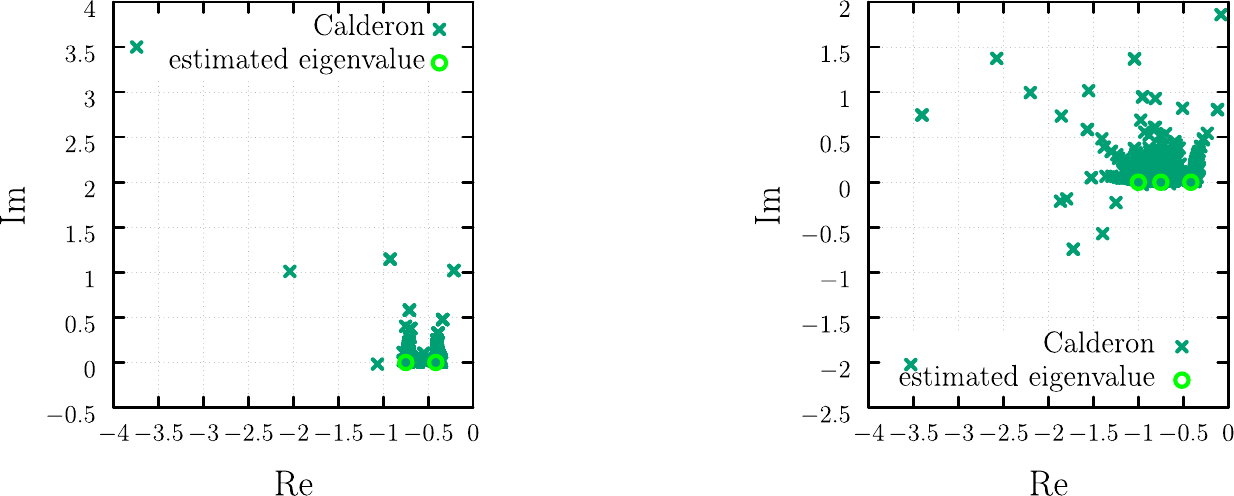}
    \caption{Eigenvalue distribution of the Calderon-preconditioned $\mathsf{A}^2$ for the configuration in Figure~\ref{fig:01kyu} (left) and that in Figure~\ref{fig:02cube} (right). Theoretically derived clustering points of the operator $\mathcal{A}^2$ from \eqref{eq:lambda} are also plotted.}
    \label{fig:0102ev}
   \end{center}
 \end{figure}

Here, we examine the eigenvalue distribution of the square of the coefficient matrix obtained by discretising the operator $\mathcal{A}$ in \eqref{eq:Ax=b} using a collocation method with constant elements. Figure~\ref{fig:0102ev} shows the distributions for the configurations in Figures~\ref{fig:01kyu} (with 3600 boundary elements) and \ref{fig:02cube} (1800 elements). For both cases, the angular frequency and material parameters are set to $\omega = 1$, $\varepsilon_1=1$, $\varepsilon_2 = 2$, and $\varepsilon_3 = 3$.
In the figure, the accumulation points of the eigenvalues of the operator ${\cal A}^2$ are also plotted. These can be evaluated using \eqref{eq:lambda} as $\tfrac{\alpha_1(1+\varepsilon_2)}{4} = -0.75$ and $\tfrac{\alpha_1(1+\varepsilon_2/\varepsilon_3)}{4} \approx -0.417$ for the case of the concentric spheres corresponding to Figure~\ref{fig:01kyu}, and $\tfrac{\alpha_1(1+\varepsilon_2)}{4} = -0.75$, $\tfrac{\alpha_1(1+\varepsilon_3)}{4} = -1$, and $\tfrac{\alpha_1(1+\varepsilon_2/\varepsilon_3)}{4} \approx -0.417$ for Figure~\ref{fig:02cube}. As observed, the eigenvalues of the coefficient matrix cluster around the predicted points, leading to faster convergence of GMRES.

\paragraph{Acknowledgements}
The point Jacobi-type preconditioning used in this paper benefited greatly from discussions with Dr~Yasuhiro Matsumoto of the Institute of Science Tokyo. The authors also acknowledge that this work was supported by JSPS KAKENHI Grant Number 21K19764. This work is also partially supported by ``Joint Usage/Research Center for Interdisciplinary Large-scale Information Infrastructures (JHPCN)'' in Japan (Project ID: jh240031 and jh250045). 

\bibliographystyle{elsarticle-num}
\bibliography{ref}

\end{document}